\documentclass[12pt,leqno]{article}
\usepackage[intlimits]{amsmath}\usepackage{amssymb}\usepackage{amsthm}

\usepackage[T2A]{fontenc}\usepackage[cp1251]{inputenc}\usepackage[russian]{babel}
\usepackage{srcltx}

\newcommand{\sign}{\mathop{\rm sign}\nolimits}
\renewcommand{\Re}{\mathop{\rm Re}\nolimits}

\renewcommand{\C}{\mathbb{C}}

\newcommand{\N}{\mathbb{N}}
\newcommand{\R}{\mathbb{R}}
\newcommand{\Z}{\mathbb{Z}}
\newcommand{\T}{\mathcal {T}}
\newcommand{\W}{{\bf W}}

\newtheorem{theorem}{\sc Теорема}[subsection]
\newtheorem{lemma}{\sc Лемма}[subsection]
\newtheorem{sledstvie}{\sc Следствие}[subsection]
\numberwithin{equation}{subsection}

\theoremstyle{remark}
\newtheorem{remark}{\sc Замечание}[subsection]
\newtheorem{example}{\sc Пример}[subsection]
\newtheorem{definition}{\sc Определение}[subsection]

\author{В. П. Заставный}
\title{Точная оценка приближения некоторых классов дифференцируемых функций сверточными операторами}

\begin{document}
\thispagestyle{empty}
\begin{center}
\textbf{\large Exact estimation of an approximation of some classes
of differentiable functions by convolution operators}
\end{center}
\bigskip

{\bf Abstract.} Nikol'skii known theorem for the kernels satisfying
a condition $A^*_n$, is proved and for kernels from wider class.
 Explicit formulas for calculating the value of an approximation
  of classes $ \W ^ {r, \beta} _ {p, n} $  by convolution operators of  special form are obtained.
 Here $ \beta\in\Z $, $r> 0 $, $n\in\N $, and $p=1 $ or $p =\infty $.
As particular cases obtained explicit formulas for value of an
approximation of the indicated classes generalized Abel-Poisson
means, biharmonic operators of Poisson, Cesaro and Riesz means.
  In some cases for value of an
approximation of the indicated classes asymptotic expansions on
parameter are found. In case of natural $r$ some results have been
obtained in  works Nikol'skii, Nagy, Timan, Telyakovskii, Baskakov,
Falaleev, Kharkevich and other mathematicians.


\medskip
  V.P. Zastavnyi\\
  Donetsk National University \\
  Universitetskaya str., 24 \\
  Donetsk 83001 Ukraine \\
  zastavn@rambler.ru

\vspace{2cm} \begin{center} \textbf{\large Точная оценка приближения
некоторых классов дифференцируемых функций сверточными операторами }
\end{center}
\bigskip

{\bf Аннотация.} Известная теорема Никольского
 для ядер, удовлетворяющих условию $A^*_n$,
 доказана и для ядер из более широкого класса.
 Получены явные  формулы для вычисления величины приближения
   классов $\W^{r,\beta}_{p,n}$ сверточными операторами    специального вида.
 Здесь   $\beta\in\Z$, $r>0$, $n\in\N$, а  $p=1$ или
 $p=\infty$.
  Как частные
 случаи получаются явные формулы для величины приближения указанных классов
 обобщенными средними Абеля-Пуассона, бигармоническими операторами Пуассона, средними  Рисса и Чезаро.
  В некоторых случаях для величины приближения указанных
 классов найдены асимптотические разложения по параметру. В случае
 натурального $r$ некоторые результаты  были получены в работах
 Никольского, Надя, Тимана, Теляковского, Баскакова, Фалалеева, Харкевича и
 других математиков.

 Ключевые слова: теорема Никольского, приближение классов функций, средние Абеля-Пуассона, бигармонические операторы Пуассона,
 средние Рисса и Чезаро,  асимптотическое
 разложение, кратно монотонная функция, функция Гурвица.

\medskip
  В.П. Заставный\\
  Донецкий Национальный Университет \\
  Университетская, 24   \\
  Донецк, 34001, Украина \\
  zastavn@rambler.ru

\newpage

\setcounter{page}{1}

\date{}
 \maketitle
\subsection{Введение}

 Пусть $L_p=L_p(-\pi,\pi)$, $1\le p\le\infty$, классы
$2\pi$-периодических вещественнозначных
 измеримых функций с конечной нормой
$||f||_p=\left(\int_{-\pi}^{\pi}|f(t)|^p dt\right)^{1/p}$ при $1\le
p<\infty$ и $||f||_{\infty}={\rm ess}\sup\{|f(t)|:t\in[-\pi,\pi]\}$.
Коэффициенты Фурье функции $\varphi\in L_1$ определяются по формуле
 $$\widehat{\varphi}(k)=\frac{1}{2\pi}\int_{-\pi}^{\pi}\varphi(t)e^{-ikt}\,dt\;,\;k\in\Z\;.$$
 Пусть $\T_n$ - множество тригонометрических полиномов вида
$$
\T_n=\left\{\frac{\alpha_0}{2}+\sum_{k=1}^{n-1}\alpha_k\cos
kt+\beta_k\sin kt: \alpha_k,\beta_k\in\R\right\}\,,\,n\in\N \;,
$$
и
 $ H_{p}^{0}=\left\{\varphi\in L_p:||\varphi||_p\le 1\right\}$,
 $ H_{p}^{n}=\left\{\varphi\in H_{p}^{0}:\widehat{\varphi}(k)=0\,,\,|k|\le n-1\,,k\in\Z\right\}$,
 $n\in\N$.
Очевидно $H_{p}^{s+1}\subset H_{p}^{s}$ при всех $s\in\Z_+$ и
$H_{p}^{n}\bot\T_n$, $n\in\N$.
\\
 По функции
$K\in L_1$ определим класс функций $\W_{p,n}(K)$:
\begin{equation}\label{W}
  \W_{p,n}(K):=\left\{f(x)=(\varphi*K)(x):=\frac{1}{2\pi}\int_{-\pi}^{\pi}\varphi(x-t)K(t)\,dt\;,\;\varphi\in
H_{p}^{n}\right\}\,,\,n\in\Z_+\,,
\end{equation}
   Очевидно $\W_{p,s+1}(K)\subset \W_{p,s}(K)$ при всех $s\in\Z_+$.
   Известно, что $\W_{p,n}(K)\subset L_p$ при всех $1\le p\le\infty$, а если $p=\infty$ или
 $K\in L_\infty$, то $\W_{p,n}(K)\subset C(\R)$ (см., например,~\cite[Гл.~4]{Korn}).

 Наилучшее приближение функции
$f$ тригонометрическими полиномами определяется по формуле
$$
 E_n(f)_p=\inf_{T\in\T_n}||f-T||_p\,,\,n\in\N \;.
 $$

  Известно, что (см., например,~\cite[(4.20), (4.24), (4.23)]{Korn})
   при $1\le p\le\infty$ справедливы соотношения
\begin{equation}\label{f1}
\sup_{f\in\W_{p,0}(K)}E_n(f)_p=\sup_{f\in\W_{p',n}(K)}||f||_{p'}\le
  \frac{1}{2\pi}\,E_n(K)_1\;,\;n\in\N\;,\;
  \frac{1}{p}+\frac{1}{p\,'}=1\,,
\end{equation}
\begin{equation}\label{f2}
\sup_{f\in\W_{1,0}(K)}E_n(f)_1=\sup_{f\in\W_{\infty,n}(K)}||f||_\infty=
\frac{1}{2\pi}\,E_n(K)_1\;,\;n\in\N\,.
\end{equation}
\begin{definition}
Говорят, что функция $K\in L_1$ удовлетворяет условию Никольского
$A_n^*$, $n\in\N$, если существуют натуральное $n_*\ge n$ и
тригонометрический полином $T^*\in\T_n$ такие, что   для функции
$\varphi_*(t)=\sign(K(t)-T^*(t))$ почти всюду\footnote{Здесь и далее
под {\it почти всюду} мы подразумеваем почти всюду относительно меры
Лебега.} выполняется равенство
$\varphi_*(t+{\pi}/{n_*})=-\varphi_*(t)$.
\end{definition}
\begin{theorem}[Никольский~(1946)~\cite{Nik1946}]\label{thNik}
Если при некотором  $n\in\N$ ядро $K\in L_1$ удовлетворяет условию
$A_n^*$ и полином $T^*\in\T_n$ из этого условия, то
 для всех $s=0,1,\ldots,n$ имеют
место соотношения
\begin{equation}\label{Nik1}
\sup_{f\in\W_{\infty,s}(K)}E_n(f)_\infty=\sup_{f\in\W_{\infty,n}(K)}||f||_\infty=
\frac{1}{2\pi}\,E_n(K)_1=\frac{1}{2\pi}\,||K-T^*||_1\;,
\end{equation}
\begin{equation}\label{Nik2}
\sup_{f\in\W_{1,s}(K)}E_n(f)_1=\sup_{f\in\W_{1,n}(K)}||f||_1=
\frac{1}{2\pi}\,E_n(K)_1=\frac{1}{2\pi}\,||K-T^*||_1\;.
\end{equation}
\end{theorem}
Теореме Никольского предшествовали исследования Колмогорова, Фавара,
Ахиезера, Крейна, Надя (более подробно см.~\cite{Nik1946}).
 \begin{definition}
Мы говорим, что функция $K\in L_1$ удовлетворяет условию $B_n^*$,
$n\in\N$, если существуют тригонометрический полином $T^*\in\T_n$,
функция $\varphi_*\in L_\infty$ и натуральное $n_*\ge n$ такие, что
почти всюду  выполняются соотношения $|\varphi_*(t)|\le 1$,
$\varphi_*(t)(K(t)-T^*(t))=|K(t)-T^*(t)|$ и
$\varphi_*(t+{\pi}/{n_*})=-\varphi_*(t)$.
 \end{definition}
 Если ядро $K$ удовлетворяет $A_n^*$ условию, то оно удовлетворяет и
$B_n^*$ условию. Обратное, вообще говоря, не верно
 (см. замечание~\ref{re0}).
   Отметим, что в условии $B_n^*$, в отличии от условия $A_n^*$,
     нам не важно на каком множестве (нулевой или положительной меры) обращается в ноль
 разность $K(t)-T^*(t)$.
 В \S~\ref{par_Nik} данной работы теорема Никольского доказана для
 ядер, которые удовлетворяют более общему условию $B^*_n$. Приведены
 как известные примеры таких ядер (ядра Надя~\cite{Nik1946,Nagy}) так и новые.

 В \S~\ref{par_Wr} рассматривается следующая задача о приближении классов
 функций сверточными операторами.
Пусть $\psi,g\in L_1$ и $1\le p\le \infty$. По функции $g\in L_1$
определим оператор
\begin{equation}\label{G}
  G(f)(x):=\frac{1}{2\pi}\int_{-\pi}^{\pi}f(x-t)g(t)\,dt\;,\;f\in \W_{p,n}(\psi)\,.
\end{equation}
 Очевидно
\begin{equation}
  f(x)-G(f)(x)=\frac{1}{2\pi}\int_{-\pi}^{\pi}\varphi(x-t)K(t)\,dt\;,\;f\in \W_{p,n}(\psi)\;,
\end{equation}
где $\varphi$ - соответствующая функция из $H_{p}^{n}$
(см.~\eqref{W}), а
\begin{equation}\label{K}
 K(x)=\psi(x)-(\psi*g)(x)\in L_1\;;\;K(x)\sim
\sum_{k}\widehat{\psi}(k)(1-\widehat{g}(k))e^{ikx}\;.
\end{equation}
Очевидно $f-G(f)\in L_p$ для любой $f\in\W_{p,n}(\psi)$, а если
$p=\infty$ или
 $\psi\in L_\infty$, то $f-G(f)\in C(\R)$. Поэтому при
$1\le q\le p\le\infty$ имеет смысл следующая величина (приближение
класса $\W_{p,n}(\psi)$ оператором $G$)
\begin{equation}
E(\W_{p,n}(\psi);G)_q:=\sup_{f\in\W_{p,n}(\psi)}||f-G(f)||_q=
\sup_{f\in\W_{p,n}(K)}||f||_q\;.
\end{equation}
Из~\eqref{f2} вытекает равенство
\begin{equation*}
E(\W_{\infty,n}(\psi);G)_\infty=\frac{1}{2\pi}\,E_n(K)_1\;,\;n\in\N\,.
\end{equation*}
Если дополнительно для ядра $K\in L_1$ выполнено условие $B_n^*$, то
(см. теорему~\ref{thBn}) справедливы равенства
\begin{equation*}
E(\W_{\infty,n}(\psi);G)_\infty= E(\W_{1,n}(\psi);G)_1=
\frac{1}{2\pi}\,E_n(K)_1=\frac{1}{2\pi}\,||K-T^*||_1\;.
\end{equation*}
 В качестве функции $\psi$ берем
\begin{equation}\label{psi}
  \psi_{r,\beta}(t)=\sum_{k\ne 0}\frac{e^{-i\beta\pi\sign
  k/2}}{|k|^r}\,e^{ikt}=
  \sum_{k=1}^{\infty}\frac{2\cos\left(kt-\frac{\beta\pi}{2}\right)}{k^r}
  \;,\;r>0\;,\;\beta\in\R\;.
\end{equation}
Отметим, что при $r=1$ справедливо равенство (см., например,
\cite[гл.~I, \S~2]{Zig})
\begin{equation}\label{psi_1}
  \psi_{1,\beta}(t)=-2\cos\frac{\beta\pi}{2}\,\ln\left(2\sin\frac{t}{2}\right)+
  \sin\frac{\beta\pi}{2}\,(\pi-t)\;,\;0<t<2\pi\;.
\end{equation}
 Известно, что  $\psi_{r,\beta}\in L_1$ (см., например,
\cite[гл.~V]{Zig} или \cite[гл.~7]{Edv}). В этом случае получаются
хорошо известные классы
$\W^{r,\beta}_{p,n}:=\W_{p,n}(\psi_{r,\beta})$. В частных случаях,
когда $n=1$, $\beta=r$ или $\beta=r+1$, получаются классы
$W^r_p:=W^{r,r}_{p,1}$ и $\widetilde{W}^r_p:=W^{r,r+1}_{p,1}$.

 В \S~\ref{par_Gad} сформулированы общие теоремы о вычислении величины $E(\W^{r,\beta}_{p,n};G)_p$ в случае,
 когда $\beta\in\Z$, $p=1$ или $p=\infty$, а оператор $G=G_{\alpha,\delta}$ в~\eqref{G}
порожден функцией $g=g_{\alpha,\delta}$, ряд Фурье которой имеет вид
\begin{equation}
g_{\alpha,\delta}(t)\sim \sum_{k\in\Z}h(|k|^\alpha\delta)e^{ikt}
\;,\;\alpha>0\,,\,\delta>0\;.
\end{equation}
Здесь $h(t)$ некоторая функция, заданная при $t\ge 0$. Если
$h(t)=e^{-t}$, то получаем операторы Абеля-Пуассона. Если
$h(t)=(1-t)^{\mu}_+$, $\mu>0$, то получаем средние
Рисса.\footnote{Здесь и далее  $t_+=t$, если $t>0$ и $t_+=0$, если
$t\le 0$.}

 В \S~\ref{par_Gadg} сформулированы общие теоремы о вычислении величины $E(\W^{r,\beta}_{p,n};G)_p$ в случае,
 когда $\beta\in\Z$, $p=1$ или $p=\infty$, а оператор $G=G_{\alpha,\delta,\gamma}$ в~\eqref{G}
порожден функцией $g=g_{\alpha,\delta,\gamma}$, ряд Фурье которой
имеет вид
\begin{equation}
g_{\alpha,\delta,\gamma}(t)\sim
\sum_{k\in\Z}(1+|k|^\alpha\delta\gamma)h(|k|^\alpha\delta)e^{ikt}
\;,\;\alpha>0\,,\,\delta>0\,,\,\gamma\in\R\;.
\end{equation}
Если $h(t)=e^{-t}$, $\alpha=1$, $\delta>0$ и
$\gamma=({1-e^{-2\delta}})/({2\delta})$, то получаем бигармонический
оператор Пуассона.

 В \S~\ref{par_Examp} для $\beta\in\Z$, $p=1$ или $p=\infty$, получены точные значения приближения классов $\W^{r,\beta}_{p,n}$
 операторами $G_{\alpha,\delta,\gamma}$, когда $h(t)=e^{-t}$, $h(t)=(1+t)^{-\mu}$, $h(t)=(1-t)_+^\mu$, $\mu>0$.
 В первых двух случаях для величин $E(\W^{r,\beta}_{p,n};G_{\alpha,\delta})_p$ найдены  асимптотические разложения по
 степеням $\delta$.

  В \S~\ref{par_Chesaro} для $\beta\in\Z$, $p=1$ или $p=\infty$,
   получены точные значения приближения классов $\W^{r,\beta}_{p,1}$
 средними Чезаро $\sigma_m^\alpha$, при $\alpha\ge 1$,
  а в \S~\ref{par_Polinom}  получены точные значения приближения этих классов
 средними типа Рисса и Чезаро.

 В \S~\ref{par_Vspom} доказаны вспомогательные утверждения.
 В \S~\ref{par_Proof} приведены доказательства теорем из  \S~\ref{par_Gad} и \S~\ref{par_Gadg}.

\subsection{Наилучшее приближение на классах сверток}\label{par_Nik}
\subsubsection{Теорема Никольского для  ядер с условием $B^*_n$}
\begin{lemma}\label{le1}
Предположим, что при некотором  $n\in\N$ для ядра $K\in L_1$
существуют тригонометрический полином $T^*\in\T_n$ и функция
$\varphi_*\in L_\infty$ такие, что $\varphi_*\bot\T_n$ и соотношения
$|\varphi_*(t)|\le 1$ и $\varphi_*(t)(K(t)-T^*(t))= |K(t)-T^*(t)|$
выполняются почти всюду на $(-\pi,\pi)$. Тогда
$E_n(K)_1=||K-T^*||_1$.
\end{lemma}
 Доказательство вытекает из следующих очевидных соотношений,
 справедливых для всех $T\in\T_n$:
 \begin{equation*}
  ||K-T^*||_1=
  \int_{-\pi}^{\pi}(K(t)-T^*(t))\varphi_*(t)\,dt=
  \int_{-\pi}^{\pi}(K(t)-T(t))\varphi_*(t)\,dt\le
  ||K-T||_1\,.
 \end{equation*}
 То, что условия в лемме~\ref{le1} являются и необходимыми для элемента
 наилучшего приближения, доказано в~\cite[Теорема 5.2.5]{TrBel}.
 \begin{remark}\label{re_le}
 Пусть  для ядра $K$
выполнено условие $B^*_n$, а полином $T^*\in\T_n$,
 функция $\varphi_*\in L_\infty$ и   натуральное $n_*\ge n$ из этого условия.
 Из неравенства $n_*\ge n$ и условия
$\varphi_*(t+{\pi}/{n_*})=-\varphi_*(t)$  вытекает,
  что функция $\varphi_*$  ортогональна многочленам из
  $\T_n$. Тогда из леммы~\ref{le1} вытекает  равенство  $E_n(K)_1=||K-T^*||_1$.
 \end{remark}
 \begin{example}\label{ex_babenko}
Пусть функция $K\in L_1$ при некоторых значениях $a,A\in\R$,
$\delta\in(0,\pi/(2n)]$, $n\in\N$, удовлетворяет условиям: $K(t)=A$
почти всюду на
  $(-\pi+a,\pi+a)\setminus(a-\delta,a+\delta)$ и
 $K(t)\ge A$ (или $K(t)\le A$)
почти всюду на $(a-\delta,a+\delta)$. Тогда $K$ очевидно
удовлетворяет условию $B^*_n$ с $n^*=n$, $T^*\equiv  A$,
   $\varphi_*(t)=\sign(\cos n(t-a))$ (или $\varphi_*(t)=-\sign(\cos n(t-a))$)
    и, значит, $E_n(K)_1=||K-T^*||_1$.

 В качестве простого примера рассмотрим функцию $K(t)=\chi_h(t)$,
 $t\in(-\pi,\pi)$, где $\chi_h$ - характеристическая функция
 интервала $(-h,h)$, $0<h<\pi$. Если $0<h\le \pi/(2n)$, $n\in\N$, то
 $E_n(K)_1=2h$ (см. пример~\ref{ex_babenko} при $a=A=0$, $\delta=h$).
 Если $0<\pi-h\le \pi/(2n)$, $n\in\N$, то
 $E_n(K)_1=2(\pi-h)$ (см. пример~\ref{ex_babenko} при $a=\pi$, $A=1$,
 $\delta=\pi-h$). Другое доказательство этих равенств, основанное на
 двойственности, содержится в работе~\cite{Babenko_Kryakin_2006} (см. также \cite[\S~5]{Babenko_Kryakin}).
 \end{example}
\begin{theorem}\label{thBn}
Предположим, что при некотором  $n\in\N$ для ядра $K\in L_1$
выполнено условие $B_n^*$, а полином $T^*\in\T_n$  из этого условия.
Тогда   для всех $s=0,1,\ldots,n$ справедливы равенства~\eqref{Nik1}
и~\eqref{Nik2}.
\end{theorem}
\begin{remark}\label{re0}
Если ядро $K$ удовлетворяет $A_n^*$ условию, то оно удовлетворяет и
$B_n^*$ условию. Обратное, вообще говоря, не верно. Условие $A_n^*$
может не выполняться в  случае, когда ядро $K(t)$ совпадает с
полиномом наилучшего приближения в $L_1$  на множестве положительной
меры. Например, функция $K(t)=(\alpha-|t|)_+$, $t\in(-\pi,\pi)$ при
$0<\alpha<{\pi}/({2n})$, $n\in\N$, очевидно удовлетворяет условию
$B_n^*$ при $T^*=0$, $n_*=n$ и $\varphi_*(t)=\sign(\cos nt)$. Если
предположить, что функция $K(t)$  удовлетворяет условию $A_n^*$ с
некоторым многочленом $T_*\in\T_n$, то в силу теоремы Джексона о
единственности многочлена наилучшего приближения в $L_1$ для
непрерывных функций (см., например, \cite{Jackson},
\cite[\S~49]{AhAppr}, \cite[\S~2.4]{TimanAF}) получим, что
$T_*=T^*$.
  Тогда функция $\sign K(t)$ должна быть ортогональна многочленам из
  $\T_n$, но $\int_{-\pi}^{\pi}\sign K(t)\,dt=2\alpha\ne 0$.
  Поэтому функция $K(t)$ не удовлетворяет условию $A_n^*$.
  Нетрудно показать, что любое непрерывное ядро, которое
  удовлетворяет условию $A_n^*$, можно исправить на множестве
  положительной меры так, чтобы исправленное ядро удовлетворяло  условию
  $B_n^*$, но не удовлетворяло условию $A_n^*$.
\end{remark}
{\sc Доказательство теоремы~\ref{thBn}.} Пусть  для ядра $K$
выполнено условие $B^*_n$, а полином $T^*\in\T_n$,
 функция $\varphi_*\in L_\infty$ и   натуральное $n_*\ge n$ из этого условия.
  В силу замечания~\ref{re_le} ядро
$K$ удовлетворяет условиям леммы~\ref{le1}. Докажем сначала
равенства
\begin{equation}\label{Nik3}
\sup_{f\in\W_{\infty,s}(K)}E_n(f)_\infty=
\sup_{f\in\W_{1,n}(K)}||f||_1=
\sup_{f\in\W_{\infty,n}(K)}||f||_\infty=
\frac{1}{2\pi}\,E_n(K)_1=\frac{1}{2\pi}\,||K-T^*||_1\;.
\end{equation}
Доказательство точно такое же как и соответствующее доказательство в
теореме Никольского в~\cite[Теорема~4.3.3]{Korn}. Так как
последовательность, стоящая в левой части~\eqref{Nik3} очевидно
убывает по $s\in\Z_+$, то в силу леммы~\ref{le1} и соотношений
\eqref{f1} при $p=\infty$ и \eqref{f2} достаточно доказать
неравенство
\begin{equation}\label{ner1}
\sup_{f\in\W_{\infty,n}(K)}E_n(f)_\infty\ge
\frac{1}{2\pi}\,E_n(K)_1\;.
\end{equation}
Если $E_n(K)_1=0$, то  неравенство очевидно. Поэтому считаем
$E_n(K)_1>0$.

Берем функцию
\begin{equation}
f_*(x):=\frac{1}{2\pi}\int_{-\pi}^{\pi}\varphi_*(t-x)K(t)\,dt=
        \frac{1}{2\pi}\int_{-\pi}^{\pi}\varphi_*(t-x)(K(t)-T^*(t))\,dt\;.
\end{equation}
  Очевидно $f_*\in\W_{\infty,n}(K)$
   с функцией $\varphi(t)=\varphi_*(-t)\in H_{\infty}^{n}$ и
 $f_*\in C(\R)$. Кроме того, функция $f_*$ имеет период
 ${2\pi}/{n_*}$ и для всех $x\in\R$ выполняются соотношения
 $$
 |f_*(x)|\le f_*(0)=
 \frac{1}{2\pi}\int_{-\pi}^{\pi}|K(t)-T^*(t)|\,dt=\frac{1}{2\pi}\,E_n(K)_1\;\;;\;\;
 f_*\left(x+{\pi}/{n_*}\right)=-f_*(x)\;.
 $$
 Поэтому $f_*\left({\pi j}/{n_*}\right)=(-1)^jf_*(0)=(-1)^j||f_*||_C$
 при всех $j\in\Z$.
 Таким образом функция $f_*$  в точках ${\pi j}/{n_*}$, $j\in\Z$,
 принимает наибольшее по абсолютной величине значения,
 последовательно меняя знак. Так как $n_*\ge n$, то этих точек на
 $[-\pi,\pi)$ не меньше $2n$ и по теореме Чебышева
 (см., например, \cite[\S~3.2]{Korn} или~\cite[\S~5.3]{TrBel})
 $$
 E_n(f_*)_\infty=E_n(f_*)_C=||f_*||_C=\frac{1}{2\pi}\,E_n(K)_1\;.
 $$
Неравенство~\eqref{ner1} доказано и, значит, доказаны
равенства~\eqref{Nik3}.

Далее точно так же как и в~\cite[Теорема~4.3.3]{Korn} можно
доказать, что (см. доказательство неравенства (4.32) из~\cite{Korn})
\begin{equation*}\label{ner2}
\sup_{f\in\W_{1,n}(K)}E_n(f)_1\ge \frac{1}{2\pi}\,E_n(K)_1\;.
\end{equation*}
Отсюда, учитывая равенство~\eqref{f2}, вытекает, что для всех
$s=0,1,\ldots,n$ справедливы равенства
\begin{equation}\label{Nik4}
\sup_{f\in\W_{1,s}(K)}E_n(f)_1=
\frac{1}{2\pi}\,E_n(K)_1=\frac{1}{2\pi}\,||K-T^*||_1\;.
\end{equation}
Теорема~\ref{thBn} доказана.
\begin{remark}
 Если две функции
$K_1,K_2\in L_1$ совпадают почти всюду и функция $K_1$ при некотором
$n\in\N$ удовлетворяет условию $B_n^*$, то и функция $K_2$,
очевидно, удовлетворяет этому условию.
\end{remark}
Исследования  наилучших приближений на классах сверток
тригонометрическими полиномами содержатся, например, в работах
Дзядыка~\cite{Dzyadyk_53,Dzyadyk}, Стечкина~\cite{Stechkin_56},
 Ефимова~\cite{Ephimov},
Сунь Юн-шеня~\cite{Sun_59}, Стечкина и
Теляковского~\cite{Stech_Tel_TrMIAN_67}, Моторного~\cite{Mornyi},
Шевалдина~\cite{Shevaldin}, Сердюка~\cite{Serdyuk},
Покровского~\cite{Pokrovskiy}.
 Отметим также  обзор Теляковского~\cite{Tel_TrMIAN_88}, работы Бабенко и
Крякина~\cite{Babenko_Kryakin_2006,Babenko_Kryakin} о приближении
характеристической функции интервала.
 Во всех указанных
работах изложена история вопроса и имеется большой список литературы
по этой тематике.

\subsubsection{Примеры  ядер с условием $B^*_n$}
Предположим, что при некотором  $n\in\N$ для ядра $K\in L_1$
выполнено условие $B_n^*$ с некоторой  функцией $\varphi_*\in
L_\infty$ и
$$
K(t)\sim \frac{\mu_0}{2}+\sum_{k=1}^{\infty}\mu_k\cos
kt+\lambda_k\sin kt= \sum_{k=-\infty}^{+\infty}c_k e^{ikt}\,,
$$
 $$
 \varphi_*(t)\sim
 \frac{\widetilde{\mu}_0}{2}+\sum_{k=1}^{\infty}\widetilde{\mu}_k\cos
kt+\widetilde{\lambda}_k\sin kt=
 \sum_{k=-\infty}^{+\infty}\widetilde{c}_k e^{ikt}\,.
 $$
 Тогда $\widetilde{\mu}_0=\widetilde{c}_0=0$. Так как свертка ${\varphi}*K$, где ${\varphi}(t)=\varphi_*(-t)$, непрерывна на $\R$, то средние
арифметические ее ряда Фурье сходятся  к ней равномерно на $\R$.
 Поэтому
 \begin{equation}\label{f}
 \begin{split}
\frac{1}{2\pi}\,E_n(K)_1 &=
\frac{1}{2\pi}\left.\int_{-\pi}^{\pi}\varphi_*(t-x)K(t)\,dt\right|_{x=0}=
 \lim_{n\to\infty}\sum_{k=-m}^{m}c_k
\widetilde{c}_{-k}\left(1-\frac{|k|}{m}\right)
\\&
=\lim_{m\to\infty}\sum_{k=1}^{m}
 \frac{\mu_k \widetilde{\mu}_{k}+\lambda_{k} \widetilde{\lambda}_{k}}{2}\left(1-\frac{k}{m}\right)
= \sum_{k=1}^{\infty}
 \frac{\mu_k \widetilde{\mu}_{k}+\lambda_{k} \widetilde{\lambda}_{k}}{2} \,.
  \end{split}
\end{equation}
Последнее равенство в~\eqref{f} справедливо,  если ряд сходится.
\begin{remark}\label{re1}
Так как $K\in L_1$, то (см. \cite[гл.~II, Теорема~8.7]{Zig})
сходится равномерно на $\R$ ряд
\begin{equation*}
\sum_{k=1}^{\infty}\frac{\mu_k\sin kx-\lambda_k\cos kx}{k}\;.
\end{equation*}
Поэтому этот ряд  является рядом Фурье своей суммы $S(x)\in C(\R)$
и, значит, при любом $n\in\N$ ряд Фурье функции
 $F(x)=\sum_{p=0}^{2n-1}(-1)^pS(x+{p\pi}/{n})$ сходится  равномерно
 на $\R$. Ряд Фурье функции $F$ легко вычисляется (аналогично как в~\cite[гл.~II, \S~1]{Zig}).
 Пусть $T_n(t)=\sum_{p=0}^{2n-1}(-1)^pe^{ip\,t}$. Тогда
\begin{equation*}
\begin{split}
F(x)&=
\sum_{k\in\Z}\widehat{S}(k)T_n\left(\frac{k\pi}{n}\right)e^{ikx}=
 2n\sum_{p\in\Z}\widehat{S}((2p+1)n)e^{i(2p+1)nx}
 \\&
 =2\sum_{p=0}^{\infty}\frac{\mu_{(2p+1)n}\sin (2p+1)nx-\lambda_{(2p+1)n}\cos (2p+1)nx}{2p+1}\;.
 \end{split}
\end{equation*}
 Полагая $x=0$ и
$x={\pi}/({2n})$, получаем сходимость  двух  рядов:
\begin{equation*}
  \sum_{k=0}^{\infty}\frac{\lambda_{(2k+1)n}}{2k+1}\;;\;
 \sum_{k=0}^{\infty}\frac{(-1)^k\mu_{(2k+1)n}}{2k+1}\;.
 \end{equation*}
\end{remark}
\begin{example}\label{ex1}
Если для ядра  $K\in L_1$
   условие $B_n^*$ выполнено  с  функцией
$$
  \varphi_*(t)=\sign(\sin nt)=
  \frac{4}{\pi}\sum_{k=0}^{\infty}\frac{\sin(2k+1)nt}{2k+1}\,,
$$
  то из теоремы~\ref{thBn}, соотношения~\eqref{f}  и  замечания~\ref{re1}
вытекает равенство
\begin{equation}\label{fs}
\sup_{f\in\W_{\infty,n}(K)}||f||_\infty=
\sup_{f\in\W_{1,n}(K)}||f||_1=
 \frac{1}{2\pi}\,E_n(K)_1=
\frac{2}{\pi}\sum_{k=0}^{\infty}\frac{\lambda_{(2k+1)n}}{2k+1}\,.
\end{equation}

 {\bf i)} Этот случай  реализуется, например, для ядер Надя~\cite{Nik1946,Nagy}
 $K\in L_1$ вида
\begin{equation}\label{nagy1}
K(t)\sim \sum_{k=1}^{\infty}\lambda_k\sin kt\,,
\end{equation}
где последовательность $\{\lambda_k\}_{k\in\N}$ убывает к нулю и
выпукла вниз (т.е. $\lambda_{k}-2\lambda_{k+1}+\lambda_{k+2}\ge 0$
при всех $k\in\N$). Так как $K\in L_1$, то  сходится ряд
$\sum_{k=1}^{\infty}{\lambda_k}/{k}$ (см. замечание~\ref{re1} при
$x=0$). В этом случае (даже без предположения выпуклости) сумма
ряда~\eqref{nagy1} $\widetilde{K}\in L_1\bigcap C(0,2\pi)$ и
ряд~\eqref{nagy1} является рядом Фурье для $\widetilde{K}$, а
частные суммы  ряда~\eqref{nagy1} сходятся в $L_1$ к функции
$\widetilde{K}$ (см., например, \cite[\S~7.3]{Edv}). Из полноты
тригонометрической системы вытекает, что $K(t)=\widetilde{K}(t)$ при
почти всех $t\in(-\pi,\pi)$. В работе~\cite[\S~2]{Nagy} (см. также
\cite[\S~7]{Nik1946}) показано, что при любом $n\in\N$ существует
нечетный полином $T^*\in\T_n$ такой, что при всех $t\in(-\pi,\pi)$
выполняется неравенство $\sin nt\,(\widetilde{K}(t)-T^*(t))\ge 0$
(приведенное в~\cite{Nagy} доказательство этого факта для
положительных $\lambda_k$ верно и для неотрицательных $\lambda_k$).
Поэтому для ядра $\widetilde{K}$, и, значит, для $K$, при любом
$n\in\N$ выполнено условие $B_n^*$ с  функцией
$\varphi_*(t)=\sign(\sin nt)$. В этом случае равенство~\eqref{fs}
справедливо при всех $n\in\N$.

 {\bf ii)} При $n=1$ этот случай  реализуется также и
для ядер $K\in L_1$ вида~\eqref{nagy1}, где $\lambda_k\ge 0$,
$k\in\N$, а последовательность $\{k\lambda_k\}_{k\in\N}$ убывает.
Тогда сходится ряд $\sum_{k=1}^{\infty}{\lambda_k}/{k}$. Поэтому
сумма ряда~\eqref{nagy1} $\widetilde{K}\in L_1\bigcap C(0,2\pi)$ и
ряд~\eqref{nagy1} является рядом Фурье для $\widetilde{K}$, а
частные суммы $\sigma_m$  ряда~\eqref{nagy1} сходятся в $L_1$ к
функции $\widetilde{K}$ и $K(t)=\widetilde{K}(t)$ при почти всех
$t\in(-\pi,\pi)$. Кроме того $\sigma_m(t)\ge 0$ при всех $m\in\N$ и
$t\in(0,\pi)$ (см., например,~\cite{Belov}).
 Отсюда следует, что $\sign(\sin t)\,\widetilde{K}(t)\ge 0$ при всех
$t\in(-\pi,\pi)$. Поэтому для ядра $\widetilde{K}$, и, значит, для
$K$,
 выполнено условие $B_1^*$ с  функцией
$\varphi_*(t)=\sign(\sin t)$ и $T^*=0$. В этом случае
равенство~\eqref{fs} справедливо при $n=1$.

 {\bf iii)} При $n=1$ этот случай  реализуется также и
для ядер $K=K_r\in L_1$ вида~\eqref{nagy1}, где
$\lambda_k={c\,(1-\nu_k)}/{k^r}$, $k\in\N$, $r=1$ или $r\ge 2$, а
$\{\nu_k\}_{k=0}^{\infty}$ такая последовательность, что  ряд
 ${\nu_0}/{2}+\sum_{k=1}^{\infty}\nu_k\cos kt$ является рядом Фурье функции $S\in L_1$ и  $S(t)\ge 0$ для
почти всех $t\in\R$ и $\nu_0\le 1$.
 Так как $E_n(cK)_1=|c|E_n(K)_1$, то для простоты рассуждений
 считаем $c=1$. Пусть
\begin{equation}\label{nagy3_r}
K_r(t)\sim \sum_{k=1}^{\infty}\frac{1-\nu_k}{k^r}\,\sin kt\,.
\end{equation}
Если $r>1$, то ряд~\eqref{nagy3_r} сходится равномерно к функции
 $\widetilde{K}_r\in  C(\R)$  и, значит, этот ряд является рядом Фурье
своей суммы. Поэтому  $K_r(t)=\widetilde{K}_r(t)$ при почти всех
$t\in(-\pi,\pi)$.

 Пусть  $r=1$.
В этом случае
\begin{equation}\label{nagy3}
 K_1(t)\sim \sum_{k=1}^{\infty}\frac{\sin kt}{k}-\sum_{k=1}^{\infty}\frac{\nu_k}{k}\,\sin kt\,.
\end{equation}
Первый ряд в~\eqref{nagy3} является рядом Фурье своей суммы, которая
при $t\in(0,2\pi)$ равна ${(\pi-t)}/{2}$. Второй ряд сходится
равномерно на $\R$ к функции
$F(t)=\int_{0}^{t}(S(x)-{\nu_0}/{2})\,dx$ (см. \cite[гл.~II,
Теоремы~2.5 и 8.7]{Zig}) и, значит, этот ряд является рядом Фурье
своей суммы $F(t)$.
  Поэтому сумма ряда~\eqref{nagy3}
$\widetilde{K}_1\in L_{\infty}\bigcap C(0,2\pi)$, а
ряд~\eqref{nagy3} является рядом Фурье для $\widetilde{K}_1$  и,
значит, $K_1(t)=\widetilde{K}_1(t)$ при почти всех $t\in(-\pi,\pi)$.
Кроме того
$$\widetilde{K}_1(t)=\frac{\pi}{2}-\frac{(1-\nu_0)t}{2}-\int_{0}^{t}S(x)\,dx\;,\;t\in(0,2\pi)\;.$$
 Поэтому функция $\widetilde{K}_1(t)$
убывает на $(0,2\pi)$ и $\widetilde{K}_1(t)\ge
\widetilde{K}_1(\pi)=0$ при всех $t\in(0,\pi)$.

Далее воспользуемся следующим  утверждением: {\it
 Пусть
$f,\varphi\in L_1$, функция $f$ нечетная и $f(t)\ge 0$ почти всюду
на $(0,\pi)$, а функция $\varphi$ четная и убывает на $(0,\pi)$.
Тогда свертка $F=f*\varphi$ является нечетной и $F(t)\ge 0$ почти
всюду на $(0,\pi)$.
 }
 Доказательство вытекает из равенства
 \begin{equation*}
 F(x)=\frac{1}{2\pi}\int_{-\pi}^{\pi}f(t)\varphi(x-t)\,dt=
 \frac{1}{2\pi}\int_{0}^{\pi}f(t)\left(\varphi(x-t)-\varphi(x+t)\right)\,dt\,,
 \end{equation*}
в котором выражение в скобках неотрицательно при всех $x\in(0,\pi)$
(при $0\le t\le\pi-x$ это вытекает из неравенств $0\le|x-t|\le
x+t\le\pi$,
 а при $\pi-x\le t\le\pi$ это следует из неравенств $0\le|x-t|=2\max\{x,t\}-(x+t)\le 2\pi-(x+t)\le\pi$).
 В этом утверждении берем $f=\widetilde{K}_1$ и $\varphi=\psi_{r-1,0}$, $r\ge 2$ (см.~\eqref{psi}).
 Очевидно $\widetilde{K}_r=\widetilde{K}_1*\psi_{r-1,0}$.
Функция $\psi_{r-1,0}$ убывает на  $(0,\pi)$. При $r=2$ это следует
из равенства~\eqref{psi_1}, а при $r>2$ из равенства
\begin{equation}\label{pr1}
   \psi'_{r-1,0}(t)=-\sum^{\infty}_{k=1}\frac{2\sin kt}{k^{r-2}}\,,\,t\in(0,2\pi)\,.
\end{equation}
 Так как
последовательность $\{k^{-r+2}\}$, $k\in\N$, убывает к нулю и
выпукла вниз, то $ \psi'_{r-1}(t)\le 0$ при
 $t\in(0,\pi)$ (см.,
например,~\cite[c.~297]{Zig} или~\cite[Гл.~IV,
Задача~6.16(г)]{Podk}).
 Таким образом, если $r=1$ или $r\ge 2$, то $\sign(\sin t)\,\widetilde{K}_r(t)\ge 0$ при всех
$t\in(-\pi,\pi)$. Поэтому для ядра $\widetilde{K}_r$, и, значит, для
$K=K_r$,
 выполнено условие $B_1^*$ с  функцией
$\varphi_*(t)=\sign(\sin t)$ и $T^*=0$. Тогда при $n=1$ справедливо
равенство~\eqref{fs}, которое в нашем случае можно записать
следующим образом
\begin{equation}\label{fs1}
\sup_{f\in\W_{\infty,1}(K)}||f||_\infty=
\sup_{f\in\W_{1,1}(K)}||f||_1=
 \frac{1}{2\pi}\,E_1(K)_1=
\frac{2|c|}{\pi}\sum_{k=0}^{\infty}\frac{1-\nu_{2k+1}}{(2k+1)^{r+1}}\,.
\end{equation}

\end{example}

 \begin{example}\label{ex2}
 Если для ядра  $K\in L_1$
   условие $B_n^*$ выполнено  с  функцией
$$
  \varphi_*(t)=\sign(\cos nt)=
  \frac{4}{\pi}\sum_{k=0}^{\infty}\frac{(-1)^k\cos(2k+1)nt}{2k+1}\,,
$$
  то из теоремы~\ref{thBn}, соотношения~\eqref{f}  и  замечания~\ref{re1}
вытекает равенство
\begin{equation}\label{fc}
\sup_{f\in\W_{\infty,n}(K)}||f||_\infty=
\sup_{f\in\W_{1,n}(K)}||f||_1=
 \frac{1}{2\pi}\,E_n(K)_1=
\frac{2}{\pi}\sum_{k=0}^{n}\frac{(-1)^k\mu_{(2k+1)n}}{2k+1}\,.
\end{equation}

 {\bf i)} Этот случай  реализуется, например, для ядер Надя~\cite{Nik1946,Nagy}
 $K\in L_1$ вида
\begin{equation}\label{nagy2}
K(t)\sim \frac{\mu_0}{2}+\sum_{k=1}^{\infty}\mu_k\cos kt\,,
\end{equation}
где последовательность $\{\mu_k\}_{k\in\N}$ убывает к нулю и при
всех $k\in\N$ выполняются неравенства
$\Delta^2\mu_k:=\mu_{k}-2\mu_{k+1}+\mu_{k+2}\ge 0$  и
$\Delta^3\mu_k:=\mu_{k}-3\mu_{k+1}+3\mu_{k+2}-\mu_{k+3}\ge 0$. В
этом случае (даже без предположения выполнения неравенств
$\Delta^3\mu_k\ge 0$) сумма ряда~\eqref{nagy2} $\widetilde{K}\in
L_1\bigcap C(0,2\pi)$ и ряд~\eqref{nagy2} является рядом Фурье для
$\widetilde{K}$ (см., например, \cite[\S~7.3]{Edv}). Из полноты
тригонометрической системы вытекает, что $K(t)=\widetilde{K}(t)$ при
почти всех $t\in(-\pi,\pi)$. В работе~\cite[\S~2]{Nagy} (см.
также~\cite[\S~7]{Nik1946}, \cite[\S~2.11.5]{TimanAF}) показано, что
при любом $n\in\N$ существует четный полином $T^*\in\T_n$ такой, что
при всех $t\in(-\pi,\pi)$ выполняется неравенство
 $\cos nt\,(\widetilde{K}(t)-T^*(t))\ge 0$
  (приведенное
в~\cite[\S~2]{Nagy} доказательство этого факта для положительных
$\mu_k$ верно и для неотрицательных $\mu_k$). Поэтому для ядра
$\widetilde{K}$, и, значит, для $K$, при любом  $n\in\N$ выполнено
условие $B_n^*$ с  функцией $\varphi_*(t)=\sign(\cos nt)$. В этом
случае равенство~\eqref{fc} справедливо при всех $n\in\N$.

 {\bf ii)} При $n=1$ этот случай  реализуется также и
для ядер $K\in L_1$ вида~\eqref{nagy2}, где $\mu_k\ge 0$, $k\in\N$,
а последовательность $\{k^2\mu_k\}_{k\in\N}$ убывает. Тогда сходится
ряд $\sum_{k=1}^{\infty}{\mu_k}$. Поэтому сумма ряда~\eqref{nagy2}
$\widetilde{K}\in C(\R)\bigcap C^1(0,2\pi)$, а ряд~\eqref{nagy2}
является рядом Фурье для $\widetilde{K}$  и, значит,
$K(t)=\widetilde{K}(t)$ при почти всех $t\in(-\pi,\pi)$. Кроме того
$\widetilde{K}'(t)\le 0$ при всех  $t\in(0,\pi)$.
 Отсюда следует, что $\sign(\cos t)\,(\widetilde{K}(t)-\widetilde{K}({\pi}/{2}))\ge 0$ при всех
$t\in(-\pi,\pi)$. Поэтому для ядра $\widetilde{K}$, и, значит, для
$K$,
 выполнено условие $B_1^*$ с  функцией
$\varphi_*(t)=\sign(\cos t)$ и $T^*=\widetilde{K}({\pi}/{2})$. В
этом случае равенство~\eqref{fc} справедливо при $n=1$.

 {\bf iii)} При $n=1$ этот случай  реализуется  и для ядер $K\in L_1$
вида~\eqref{nagy2}, где  последовательность $\{k\mu_k\}_{k\in\N}$
убывает к нулю и выпукла вниз. В этом случае и последовательность
$\{\mu_k\}_{k\in\N}$ убывает к нулю и также выпукла вниз. Поэтому
сумма ряда~\eqref{nagy2} $\widetilde{K}\in L_1\bigcap C(0,2\pi)$ и
ряд~\eqref{nagy2} является рядом Фурье для $\widetilde{K}$   и,
значит, $K(t)=\widetilde{K}(t)$ при почти всех $t\in(-\pi,\pi)$.
Кроме того, $\widetilde{K}\in C^1(0,2\pi)$ и $\widetilde{K}'(t)\le
0$ при всех  $t\in(0,\pi)$ (см., например,~\cite[c.~297]{Zig}
или~\cite[Гл.~IV, Задача~6.16(г)]{Podk}).
 Отсюда следует, что $\sign(\cos t)\,(\widetilde{K}(t)-\widetilde{K}({\pi}/{2}))\ge 0$ при всех
$t\in(-\pi,\pi)$, $t\ne 0$. Поэтому для ядра $\widetilde{K}$, и,
значит, для $K$, выполнено условие $B_1^*$ с  функцией
$\varphi_*(t)=\sign(\cos t)$ и $T^*=\widetilde{K}({\pi}/{2})$. В
этом случае равенство~\eqref{fc} справедливо при $n=1$.

 {\bf iv)} При $n=1$ этот случай  реализуется также и
для ядер $K=M_r\in L_1$ вида~\eqref{nagy2}, где
$\mu_k={c(1-\nu_k)}/{k^r}$, $k\in\N$, $r=2$ или $r\ge 3$, а
$\{\nu_k\}_{k=0}^{\infty}$ такая последовательность, что  ряд
 ${\nu_0}/{2}+\sum_{k=1}^{\infty}\nu_k\cos kt$ является рядом Фурье функции $S\in L_1$ и  $S(t)\ge 0$ для
почти всех $t\in\R$ и $\nu_0\le 1$.
 Так как $E_n(cK)_1=|c|E_n(K)_1$, то для простоты рассуждений
 считаем $c=1$. Пусть
\begin{equation}\label{nagy3_Mr}
M_r(t)\sim \sum_{k=1}^{\infty}\frac{1-\nu_k}{k^r}\,\cos kt\,.
\end{equation}
Если $r>1$, то ряд~\eqref{nagy3_Mr} сходится равномерно к функции
 $\widetilde{M}_r\in  C(\R)$  и, значит, этот ряд является рядом Фурье
своей суммы. Поэтому  $M_r(t)=\widetilde{M}_r(t)$ при почти всех
$t\in(-\pi,\pi)$.
 Если  $r=2$ или $r\ge 3$, то $\widetilde{M}'_r(t)=-\widetilde{K}_{r-1}(t)\le 0$
 при  всех $t\in(0,\pi)$ (см. пример \ref{ex1}(iii)).
 Отсюда следует, что $\sign(\cos t)\,(\widetilde{M}_r(t)-\widetilde{M}_r({\pi}/{2}))\ge 0$ при всех
$t\in(-\pi,\pi)$. Поэтому для ядра $\widetilde{M}_r$, и, значит, для
$K=M_r$, выполнено условие $B_1^*$ с  функцией
$\varphi_*(t)=\sign(\cos t)$ и $T^*=\widetilde{M}_r({\pi}/{2})$.
 Тогда при $n=1$ справедливо равенство~\eqref{fc}, которое в нашем
случае можно записать следующим образом
\begin{equation}\label{fc1}
\sup_{f\in\W_{\infty,1}(K)}||f||_\infty=
\sup_{f\in\W_{1,1}(K)}||f||_1=
 \frac{1}{2\pi}\,E_1(K)_1=
\frac{2|c|}{\pi}\sum_{k=0}^{\infty}\frac{(-1)^k(1-\nu_{2k+1})}{(2k+1)^{r+1}}\,.
\end{equation}
\end{example}
\begin{example}
{\bf i)} Пусть функция $K\in L_1$ при некотором значении $A\in\R$
удовлетворяет условиям:
 $K(t)\ge A$ почти всюду на $(0,\pi)$ и
  $K(t)\le A$ почти всюду на $(-\pi,0)$.
  Тогда $K$ очевидно удовлетворяет условию $B^*_1$ с $n^*=1$, $T^*\equiv  A$,
   $\varphi_*(t)=\sign(\sin t)$. В этом случае равенство~\eqref{fs} справедливо
   при $n=1$.

{\bf ii)}  Пусть функция $K\in L_1$ при некотором значении $A\in\R$
удовлетворяет условиям:
 $K(t)\ge A$ почти всюду на $(-{\pi}/{2},{\pi}/{2})$ и
  $K(t)\le A$ почти всюду на
  $(-\pi,\pi)\setminus(-{\pi}/{2},{\pi}/{2})$.
  Тогда $K$ очевидно удовлетворяет условию $B^*_1$ с $n^*=1$, $T^*\equiv  A$,
   $\varphi_*(t)=\sign(\cos t)$. В этом случае равенство~\eqref{fc} справедливо
   при $n=1$.
\end{example}

\begin{remark}
Отметим, что пример \ref{ex1}(iii) при нечетных $r\in\N$, а  пример
\ref{ex2}(iv) при четных $r\in\N$ были  получены другим методом в
\cite{Pych1967}.

 Приведём простые достаточные условия неотрицательности
функции $S(t)$ из  примеров \ref{ex1}(iii) и \ref{ex2}(iv). Функция
$f:\R \to \C$ называется положительно определённой на $\R$ (см.,
например, \cite[\S~6.2]{TrBel}, \cite{Ah}), если для любых $ n\in
\N$, $ \{x_k\}_{k=1}^n \subset \R$ и $ \{c_k\}_{k=1}^n \subset \C$
 выполняется неравенство $ \sum_{k,j=1}^n c_k \bar{c}_j f(x_k-x_j) \ge 0 $.
 Для таких функций   $|f(x)|\le f(0)$, $x\in\R$ и непрерывность в
нуле эквивалентна непрерывности на~$\R$. По теореме Бохнера-Хинчина
функция $f$ является положительно определённой и непрерывной на~$\R$
тогда и только тогда, когда
 $f(x)=\int_{-\infty}^{+\infty}
e^{-iux} \,d\mu (u),$
 где $\mu$ неотрицательная, конечная,
борелевская мера на $\R$. Если $f\in C(\R) \cap L(\R),$ то
положительная определенность функции $f$ эквивалентна
неотрицательности её преобразования Фурье, т.е.
$\widehat{f}(x):=\int_{-\infty}^{+\infty} f(u)e^{-iux}\,du\ge 0 $,
$x\in \R $ и в этом случае $\widehat{f}\in L(\R)$ (см. \cite[гл. I,
\S 1, следствие 1.26]{StW}).
 \begin{lemma}\label{le0}
 Пусть функция $f$ является положительно определенной и непрерывной
 на $\R$. Если ряд $\sum_{k\in\Z}f(k)e^{ikt}$ является рядом Фурье
 функции $S\in L_1$, то $S(t)\ge 0$ при почти всех $t\in\R$.
 \end{lemma}
 {\sc Доказательство.}
  Очевидно
 функция
  $\varphi_n(x):=f(x)\left(1-{|x|}/{n}\right)_+$ при любом $n\in\N$
 имеет компактный носитель и
 является положительно определенной и непрерывной на $\R$. Тогда
  $\sigma_n(S)(t)=\sum_{k\in\Z}\varphi_n(k)e^{ikt}\ge 0$ при $t\in\R$
   (см., например, \cite[\S~6]{Zast2000}). Так как средние
   арифметические   ряда Фурье функции  $S\in L_1$ сходятся в $L_1$ к функции $S$,
   то для некоторой подпоследовательности
   $\sigma_{n_k}(S)(t)\to S(t)$ при почти всех $t\in[0,2\pi]$
    и, значит,  $S(t)\ge 0$ почти всюду на $\R$. Лемма~\ref{le0}
    доказана.
\end{remark}

 \subsection{Приближение классов $\W^{r,\beta}_{p,n}$ операторами  специального вида}\label{par_Wr}
 \subsubsection{Случай  операторов $G_{\alpha,\delta}$}\label{par_Gad}

\begin{theorem}\label{thZast1}
Пусть для функции $h$ выполнены следующие два условия:\\
{\bf 1.} Функция $\displaystyle \lambda(x)={(1-h(x))}/{x}$  выпукла
вниз на $(0,+\infty)$.\\
{\bf 2.} При некоторых $\alpha\in(0,1]$ и $\delta>0$ ряд
 $
\sum_{k\in\Z}h(|k|^\alpha\delta)e^{ikt} $
  является рядом Фурье некоторой функции $g_{\alpha,\delta}\in L_1$.
\\
Пусть  $G_{\alpha,\delta}$ оператор, порожденный функцией
$g_{\alpha,\delta}\in L_1$ по формуле~\eqref{G} и $p=1$ или
$p=\infty$. Тогда имеют место следующие утверждения:\\
{\bf 1.}
 Если $\beta+1\in2\Z$, то   при любых $n\in\N$ и $r\ge\alpha$   справедливы
равенства
\begin{equation}\label{Zast1}
E(\W^{r,\beta}_{p,n};G_{\alpha,\delta})_p=
  \displaystyle\frac{4}{\pi n^r}\sum\limits_{k=0}^{\infty}\frac{1-h((2k+1)^{\alpha}n^{\alpha}\delta)}{(2k+1)^{r+1}}\;.
\end{equation}
{\bf 2.} Если $\beta\in 2\Z$, то   для $n=1$ и любых $r\ge\alpha+1$
справедливы равенства
\begin{equation}\label{Zast1_2}
E(\W^{r,\beta}_{p,n};G_{\alpha,\delta})_p= \displaystyle
\frac{4}{\pi n^r}\sum\limits_{k=0}^{\infty}\frac{1-h((2k+1)^{\alpha}n^{\alpha}\delta)}{(2k+1)^{r+1}}\,(-1)^k\;.
\end{equation}
Если дополнительно $\lambda(x)\in C^1(0,+\infty)$ и функция
$-\lambda'(x)$ выпукла вниз на $(0,+\infty)$, то
равенства~\eqref{Zast1_2} справедливы при любых $n\in\N$ и
$r\ge\alpha$.
\end{theorem}
\begin{definition}\label{defM_m}
Обозначим через $M_m$, $m\in\N$, класс функций $h\in
C^{m-1}(0,+\infty)$, для которых  функция $(-1)^{m-1}h^{(m-1)}$
неотрицательна, убывает, выпукла вниз на интервале $(0,+\infty)$ и
существует конечный предел $h(+\infty)\ge 0$.
\end{definition}
\begin{remark}\label{re_M_m}
Примеры функций, которые удовлетворяют условиям
теоремы~\ref{thZast1}, приведены в следствии~\ref{sl}, из которого
вытекает, что если
 $h\in M_{m+1}$ при некотором  $m\in\N$ и существует конечный
предел $h(+0)\le 1$, то $\lambda(x)={(1-h(x))}/{x}\in M_k$ при всех
$k=1,\ldots,m$.
\end{remark}
\begin{definition}\label{def_m(h)}
Для функции $h$, которая удовлетворяет неравенству $h(x)\le 1$ при
всех $x>0$, определим величину $m(h)$ как точную нижнюю грань тех
$\gamma\in\R$, для которых функция
 ${(1-h(x))}/{x^{\gamma}}$ убывает на $(0,+\infty)$
 Если такие $\gamma\in\R$ не существуют, то считаем, что $m(h):=+\infty$.
\end{definition}
\begin{theorem}\label{thZast}
Пусть для функции $h(x)$, $x\ge 0$, выполнены следующие два
условия:\\
{\bf 1.} $h(x)\le 1$ при  $x>0$ и $m(h)<+\infty$.
\\
{\bf 2.} При некоторых $\alpha>0$ и $\delta>0$ ряд
 $
\sum_{k\in\Z}h(|k|^\alpha\delta)e^{ikt}
$
  является рядом Фурье некоторой функции $g_{\alpha,\delta}\in L_1$.
\\
  Тогда $m(h)\ge 0$ и $m(h)>0$, если $h(x_0)>0$ в некоторой точке
  $x_0>0$.
 Пусть $G_{\alpha,\delta}$ оператор, порожденный функцией
$g_{\alpha,\delta}\in L_1$ по формуле~\eqref{G} и  $p=1$ или
$p=\infty$. Тогда имеют место следующие утверждения:\\
{\bf 1.} Если $\beta+1\in2\Z$, то равенство~\eqref{Zast1}
справедливо для $n=1$ и  любых
 $r\ge\alpha m(h)+1$.
\\
{\bf 2.} Если $\beta\in2\Z$, то равенство~\eqref{Zast1_2}
справедливо для $n=1$ и  любых
 $r\ge\alpha m(h)+2$.
\end{theorem}
\begin{remark}\label{re_m(h)}
 Если, например, функция $h$ выпукла вниз на
$(0,+\infty)$, $h(x)\le 1$ при  $x>0$ и $h(x)\not\equiv const$ на
$(0,+\infty)$,
 то из леммы~\ref{le2} вытекает, что $m(h)\in(0,1]$.
 Если дополнительно $h'$ непрерывна в точке $0$ справа, $h(0)=1$ и $h'(0)\ne 0$, то $m(h)=1$.
\end{remark}

\subsubsection{Случай  операторов $G_{\alpha,\delta,\gamma}$}\label{par_Gadg}

\begin{definition}\label{def_gamma_m}
Для функции $h\in M_{m+1}$, $m\in\N$, с условиями $0<h(+0)\le 1$ и
$h(+\infty)=0$ определим величину $\gamma_{m}(\rho,h)$, $\rho\ge 1$,
как точную верхнюю грань тех $\gamma\in\R$, для которых функция
$\lambda_{\rho,\gamma}(x)\in M_{m}$, где
\begin{equation}\label{lambda_rho}
\lambda_{\rho,\gamma}(x)=\frac{1-(1+\gamma x)h(x)}{x^{\rho}}=
  \frac{1-h(x)}{x^{\rho}}-\gamma \frac{h(x)}{x^{\rho-1}}\,.
\end{equation}
\end{definition}
\begin{remark}
Из леммы~\ref{le5} вытекает, что для любой функции $h\in M_{m+1}$,
$m\in\N$, с условиями $0<h(+0)\le 1$ и $h(+\infty)=0$, неравенство
$0\le\gamma_m(\rho,h)<+\infty$ выполняется при любом $\rho\ge 1$.
\end{remark}
\begin{remark}\label{re_gamma_m}
В работе~\cite{Zast2009_IPMM} величина $\gamma_{m}(\rho,h)$ найдена
в следующих случаях:\\
1) Если $h(t)=e^{-t}$ и $m\in\N$, то $\gamma_m(1,h)=\frac{1}{m+2}$,
$\gamma_m(\rho,h)=1$ при $\rho\ge 2$.\\
2) Если $h_\mu(t)=(1+t)^{-\mu}$, $\mu\ge 1$ и $m\in\N$,
 то $\gamma_m(\rho,h_\mu)=\mu$ при $\rho\ge 2$, $\gamma_m(1,h_\mu)=\frac{\mu+m+1}{m+2}$,
  $\gamma_m(\rho,h_1)=1$ при $\rho\ge 1$.\\
 3) Если $H_\mu(t)=(1-t)_+^{\mu}$, $\mu\ge m+1$, $m\in\N$, то
$\gamma_m(1,H_\mu)=\frac{\mu-m-1}{m+2}$ и $\gamma_m(\rho,H_\mu)=\mu$
  при $\rho\ge 2$.
\end{remark}
\begin{theorem}\label{thZast3}
Пусть для функции $h(x)$, $x\ge 0$, выполнены следующие два
условия:\\
{\bf 1.} Функция $h\in M_2$
 и $0<h(+0)\le 1$,  $h(+\infty)=0$.
 \\
 {\bf 2.}
 При некоторых $\rho\ge 1$, $\gamma\le\gamma_1(\rho,h)$, $\alpha>0$, $\delta>0$ ряд
 $
\sum_{k\in\Z}(1+|k|^\alpha\delta\gamma)h(|k|^\alpha\delta)e^{ikt} $
  является рядом Фурье некоторой функции $g_{\alpha,\delta,\gamma}\in L_1$.
\\
Пусть  $G_{\alpha,\delta,\gamma}$ оператор, порожденный функцией
$g_{\alpha,\delta,\gamma}\in L_1$ по формуле~\eqref{G} и $p=1$ или
$p=\infty$. Тогда имеют место следующие утверждения:\\
{\bf 1.}
 Если $\beta+1\in2\Z$ и выполнено одно из двух условий:
  {\bf i)} $n\in\N$, $\alpha\in(0,1]$,  $r\ge\alpha\rho$
 или
  {\bf ii)} $n=1$, $\alpha>1$,  $r\ge\alpha\rho+1$,
 то  справедливы равенства
\begin{equation}\label{Zast2_1}
E(\W^{r,\beta}_{p,n};G_{\alpha,\delta,\gamma})_p=
  \displaystyle\frac{4}{\pi n^r}\sum\limits_{k=0}^{\infty}\frac{1-(1+(2k+1)^{\alpha}n^{\alpha}\delta\gamma)h((2k+1)^{\alpha}n^{\alpha}\delta)}{(2k+1)^{r+1}}\;.
\end{equation}
{\bf 2.} Если $\beta\in 2\Z$ и выполнено одно из двух условий:
  {\bf i)} $n=1$, $\alpha\in(0,1]$,  $r\ge\alpha\rho+1$
 или
  {\bf ii)}~$n=1$, $\alpha>1$,  $r\ge\alpha\rho+2$,
 то  справедливы равенства
\begin{equation}\label{Zast2_2}
E(\W^{r,\beta}_{p,n};G_{\alpha,\delta,\gamma})_p= \displaystyle
\frac{4}{\pi n^r}\sum\limits_{k=0}^{\infty}\frac{1-(1+(2k+1)^{\alpha}n^{\alpha}\delta\gamma)h((2k+1)^{\alpha}n^{\alpha}\delta)}{(2k+1)^{r+1}}\,(-1)^k\;.
\end{equation}
Если дополнительно $h\in M_3$,  $\gamma\le\gamma_2(\rho,h)$ и
$\alpha\in(0,1]$, то равенства~\eqref{Zast2_2} справедливы при любых
$n\in\N$ и $r\ge\alpha\rho$.
\end{theorem}
В следующей теореме рассмотрен случай положительных ядер, которые
порождаются положительно определенными функциями.
\begin{theorem}\label{thZast4}
Пусть для функции $h(x)$, $x\ge 0$ выполнены следующие условия:\\
{\bf 1.} При некоторых  $\gamma\in\R$, $\alpha>0$, $\delta>0$ ряд
 $
\sum_{k\in\Z}(1+|k|^\alpha\delta\gamma)h(|k|^\alpha\delta)e^{ikt} $
  является рядом Фурье некоторой функции $g_{\alpha,\delta,\gamma}\in L_1$.
\\
{\bf 2.} Функция $(1+\gamma |t|^\alpha)h(|t|^\alpha)$ является
положительно определенной и непрерывной на $\R$ и $h(0)\le 1$.
\\
Пусть $p=1$ или $p=\infty$, а  $G_{\alpha,\delta,\gamma}$ -
оператор, порожденный функцией $g_{\alpha,\delta,\gamma}$ по
формуле~\eqref{G}. Тогда равенство~\eqref{Zast2_1} справедливо при
$n=1$, $r=1$  и $r\ge 2$, а  равенство~\eqref{Zast2_2} справедливо
при $n=1$, $r=2$  и $r\ge 3$.
\end{theorem}

\subsubsection{Примеры операторов $G_{\alpha,\delta}$ и
$G_{\alpha,\delta,\gamma}$}\label{par_Examp}
\begin{example}\label{Abel}
Пусть $h(t)=e^{-t}$.  В этом случае $m(h)=1$ (см.
замечание~\ref{re_m(h)}).
 Второе условие в
теоремах~\ref{thZast1} и \ref{thZast} очевидно выполнено для любых
$\alpha>0$ и $\delta>0$. Очевидно $h\in M_m$ при любом $m\in\N$.
Учитывая замечание~\ref{re_M_m}, получаем следующие результаты. Если
$\alpha\in(0,1]$, то равенства~\eqref{Zast1} и ~\eqref{Zast1_2}
выполняются при любых $n\in\N$, $r\ge\alpha$ и $\delta>0$. Если
$\alpha>1$, то равенство~\eqref{Zast1} выполняется при  $n=1$ и
любых $r\ge\alpha+1$ и $\delta>0$, а равенство~\eqref{Zast1_2}
выполняется при  $n=1$ и любых $r\ge\alpha+2$ и $\delta>0$.
Соответствующие операторы $G_{\alpha,\delta}$ называются обобщенными
операторами Абеля-Пуассона. Для операторов Абеля-Пуассона
($\alpha=1$) результат был известен только при $n=1$, $\beta=r\in\N$
(см.~\cite{Timan}) и  при $n=1$, $\beta-1=r\in\N$  (см.
\cite{Nagy1950} при $r=1$  и \cite{Harkevich2002_1} при $r\ge 2$).
Оба эти случая вытекают из теоремы~\ref{thZast1}.

Второе условие в теореме~\ref{thZast3}  выполнено для любых
$\alpha>0$, $\delta>0$ и $\gamma\in\R$, а значения
$\gamma_m(\rho,h)$ найдены в работе автора~\cite{Zast2009_IPMM}
 (см. замечание~\ref{re_gamma_m}).
 Поэтому теорема~\ref{thZast3}  справедлива, например, в следующих
 случаях:
 {\bf 1)} $\rho=1$, $\gamma_1(1,h)=\frac 13$, $\gamma_2(1,h)=\frac 14$.
 {\bf 2)}~$\rho=2$, $\gamma_1(2,h)=\gamma_2(2,h)=1$.

 Так как функция $e^{-|t|^\alpha}$ является положительно
 определенной на $\R$
$\iff$ $0<\alpha\le 2$, то операторы $G_{\alpha,\delta,\gamma}$ при
$\gamma=0$ будут положительными при любых $0<\alpha\le 2$ и
$\delta>0$. И наоборот, если при некотором $\alpha>0$ операторы
$G_{\alpha,\delta,0}$ будут положительными при любых  $\delta>0$, то
$0<\alpha\le 2$  (более подробно см., например, \cite{Zast2000}).
 Нетрудно показать, что функция $(1+\gamma|x|)e^{-|x|}$ является
положительно определенной на $\R$ $\iff\gamma\in[-1,1]$
(преобразование Фурье этой функции равно
${2(1+\gamma+(1-\gamma)t^2)}/{(1+t^2)^2}$).
 Поэтому операторы
$G_{1,\delta,\gamma}$ будут положительными при любых  $\delta>0$
$\iff\gamma\in[-1,1]$.
 В силу теоремы~\ref{thZast4} равенство~\eqref{Zast2_1} справедливо
 в следующих  случаях:
 {\bf 1)} $n=1$, $r=1$ или $r\ge 2$, $\gamma=0$, $0<\alpha\le 2$.
 {\bf 2)}~$n=1$, $r=1$  или $r\ge 2$, $\gamma\in[-1,1]$, $\alpha=1$.
 Отметим, что  случай {\bf 1)} при $r=1$ хорошо известен (см., например,
 \cite{Bausov1961,Baskakov,Falaleev2001}).
 Равенство~\eqref{Zast2_2} справедливо
 в следующих  случаях:
 {\bf 1)}~$n=1$, $r=2$ или $r\ge 3$, $\gamma=0$, $0<\alpha\le 2$.
 {\bf 2)}~$n=1$, $r=2$  или $r\ge 3$, $\gamma\in[-1,1]$, $\alpha=1$.

Для бигармонических операторов Пуассона ( $\alpha=1$, $\delta>0$ и
$\gamma={(1-e^{-2\delta})}/{(2\delta)}\in(0,1)$) результат был
известен только при $n=1$ и  $r-1\in\N$ (см.~\cite{Harkevich2002}),
который вытекает из теоремы~\ref{thZast3} и при $n=r=\beta=1$
(см.~\cite{Falaleev, Harkevich2000}), который вытекает из
теоремы~\ref{thZast4}.

 В случае операторов Абеля-Пуассона ($\alpha=1$)
  поиску полного асимптотического
представления рядов~\eqref{Zast1} и~\eqref{Zast1_2}  при $n=1$,
$\beta,r\in\N$ были посвящены
работы~\cite{Harkevich2002_1,Baskakov,Maley,Stark}. Полное решение
этой задачи было получено в работе автора~\cite{Zast2009_mz}.
 В случае $\alpha=2$, $n=1$, $\beta=r=1$ отметим работу
 \cite{Baskakov}.
 В общем случае  асимптотические разложения  рядов~\eqref{Zast1},
\eqref{Zast1_2}, \eqref{Zast2_1} и~\eqref{Zast2_2} с явными
коэффициентами легко вытекают из результатов работы
автора~\cite{Zast2009}. Выпишем эти разложения при $\gamma=0$ (в
случае $\gamma\ne 0$ они слишком громоздкие).

Пусть $r>0$,  $\alpha>0$. Тогда следующие асимптотические разложения
справедливы
 соответственно в случаях, когда
${r}/{\alpha}\not\in\N$ и ${r}/{\alpha}=p\in\N$:
\begin{equation}\label{ass1}
\begin{split}
&\sum_{k=0}^{\infty}\frac{1-e^{-\delta(2k+1)^\alpha}}{(2k+1)^{r+1}}
\underset{\delta\to+0}{\sim}
\\&
 -2^{-r-1}\left[
\frac{1}{\alpha}\;\Gamma\left(-{r}/{\alpha}\right)\,
      (2^\alpha\delta)^{\frac{r}{\alpha}}+
      \sum_{k=1}^{\infty}\frac{(-1)^k}{k!}\,\zeta\left(-\alpha
      k+r+1, 1/2\right)\,(2^\alpha\delta)^k\right]\;,
      \end{split}
      \end{equation}
      \begin{equation}\label{ass2}
\begin{split}
\sum_{k=0}^{\infty}\frac{1-e^{-\delta(2k+1)^\alpha}}{(2k+1)^{r+1}}
\underset{\delta\to+0}{\sim}
  &\;
   -2^{-r-1}\left[
  \frac{(-1)^p(2^\alpha\delta)^p}{\Gamma(p+1)}\left(-\frac{\ln (2^\alpha\delta)}{\alpha}+\frac{\Gamma'(p+1)}{\alpha\Gamma(p+1)}
  -\frac{\Gamma'\left( 1/2\right)}{\Gamma\left(1/2\right)}\right)\,
      \right.
      \\&
      +\left.\sum_{k=1,k\ne p}^{\infty}\frac{(-1)^k}{k!}\,\zeta\left(-\alpha
      k+r+1, 1/2\right)\,(2^\alpha\delta)^k\right]\;.
\end{split}
\end{equation}
 Если
$0<\alpha<1$, то в~\eqref{ass1} и~\eqref{ass2} имеет место знак
равенства при всех  $\delta>0$. Если $\alpha=1$, то в~\eqref{ass1}
и~\eqref{ass2} имеет место знак равенства при всех
$\delta\in(0,\pi)$. Здесь  ${\zeta}(s,a)$ - функция Гурвица с
параметром $a>0$, равная  при $\Re s>1$ сумме
$\sum_{k=0}^{\infty}(k+a)^{-s}$. Эта функция аналитически
продолжается в $\C\setminus\{1\}$.

Пусть $r+1>0$, $\alpha>0$. Тогда имеет место следующее
асимптотическое разложение
\begin{equation}\label{ass3}
\sum_{k=0}^{\infty}\frac{1-e^{-\delta(2k+1)^\alpha}}{(2k+1)^{r+1}}(-1)^k
\underset{\delta\to+0}{\sim}
     -2^{-r-1}
      \sum_{k=1}^{\infty}\frac{(-1)^k}{k!}\,\widetilde{\zeta}\left(-\alpha k+r+1, 1/2\right)\,(2^\alpha\delta)^k\;.
      \end{equation}
 Если $0<\alpha<1$, то
в~\eqref{ass3}  имеет место знак равенства при всех $\delta>0$. Если
$\alpha=1$, то в~\eqref{ass3}  имеет место знак равенства при всех
$\delta\in\left(0,{\pi}/{2}\right)$. Здесь $\widetilde{\zeta}(s,a)$
- целая функция по $s\in\C$, равная  при $\Re s>0$ сумме
$\sum_{k=0}^{\infty}(-1)^k(k+a)^{-s}$,  $a>0$.

\end{example}
\begin{example}\label{New}
Пусть $h(t)=(t+1)^{-\mu}$, $\mu>0$.   В этом случае $m(h)=1$ (см.
замечание~\ref{re_m(h)}). Второе условие в теоремах~\ref{thZast1} и
\ref{thZast} выполнено для любых $\alpha>0$ и $\delta>0$. Это
вытекает из того, что функция $h(t^\alpha\delta)$ убывает к нулю и
выпукла вниз на $(t_{\alpha,\mu,\delta},+\infty)$ при некотором
$t_{\alpha,\mu,\delta}>0$.  Очевидно $h\in M_m$ при любом $m\in\N$.
Учитывая замечание~\ref{re_M_m}, получаем следующие результаты.
Если $\alpha\in(0,1]$, то равенства~\eqref{Zast1} и ~\eqref{Zast1_2}
выполняются при любых $n\in\N$, $r\ge\alpha$ и $\delta>0$. Если
$\alpha>1$, то равенство~\eqref{Zast1} выполняется при  $n=1$ и
любых $r\ge\alpha+1$ и $\delta>0$, а равенство~\eqref{Zast1_2}
выполняется при  $n=1$ и любых $r\ge\alpha+2$ и $\delta>0$.

Второе условие в теореме~\ref{thZast3} выполнено для любых $\mu>1$,
$\alpha>0$, $\delta>0$ и $\gamma\in\R$. Это вытекает из того, что
функция
 $t^\alpha h(t^\alpha\delta)$ убывает к нулю и выпукла вниз на
$(t_{\alpha,\mu,\delta},+\infty)$ при некотором
$t_{\alpha,\mu,\delta}>0$.
 Значения
$\gamma_m(\rho,h)$ найдены в  работе~\cite{Zast2009_IPMM} (см.
замечание~\ref{re_gamma_m}).
 Поэтому теорема~\ref{thZast3}  справедлива в следующих
 случаях:
 {\bf 1)}~$\rho=1$, $\gamma_1(1,h)={(\mu+2)}/{3}$, $\gamma_2(1,h)={(\mu+3)}/{4}$, $\mu>1$.
 {\bf 2)}~$\rho=2$, $\gamma_1(2,h)=\gamma_2(2,h)=\mu$, $\mu>1$.

     Функция $(|t|^\alpha+1)^{-\mu}$, $\alpha>0$, $\mu>0$, является положительно
 определенной на $\R$
$\iff$ $0<\alpha\le 2$. Достаточность вытекает из того, что функция
$(t+1)^{-\mu}$, $\mu>0$, является вполне монотонной на
$(0,+\infty)$. Необходимость вытекает из того, что среди
положительно определенных функций только постоянная функция имеет в
нуле нулевую производную второго порядка. Поэтому операторы
$H_{\alpha,\mu,\delta,\gamma}=G$ при $\gamma=0$ будут положительными
при любых $0<\alpha\le 2$, $\mu>0$ и $\delta>0$. И наоборот, если
при некоторых $\alpha>0$, $\mu>0$ операторы
$H_{\alpha,\mu,\delta,0}$ будут положительными при любых $\delta>0$,
то $0<\alpha\le 2$ (более подробно см., например, \cite{Zast2000}).
 В силу теоремы~\ref{thZast4} равенство~\eqref{Zast2_1} справедливо
 в следующем  случае:  $n=1$, $r=1$ или $r\ge 2$, $\gamma=0$, $0<\alpha\le 2$, $\mu>0$.
Равенство~\eqref{Zast2_2} справедливо
 в следующем  случае:  $n=1$, $r=2$ или $r\ge 3$, $\gamma=0$, $0<\alpha\le 2$, $\mu>0$.

Асимптотические разложения  рядов~\eqref{Zast1}, \eqref{Zast1_2},
\eqref{Zast2_1} и~\eqref{Zast2_2} с явными коэффициентами легко
вытекают из результатов работы~\cite{Zast2009}. Выпишем эти
разложения при $\gamma=0$ (в случае $\gamma\ne 0$ они слишком
громоздкие).

  Пусть $r>0$,
$\alpha>0$ и $\mu>0$. Тогда следующие асимптотические разложения
справедливы соответственно в случаях, когда ${r}/{\alpha}\not\in\N$
и ${r}/{\alpha}=p\in\N$
\begin{equation}\label{assg1}
\begin{split}
&\sum_{k=0}^{\infty}\frac{1-((2k+1)^\alpha\delta+1)^{-\mu}}{(2k+1)^{r+1}}
\underset{\delta\to+0}{\sim}
  \\
  &
  -2^{-r-1}\left[
  \frac{\Gamma\left(-\frac{r}{\alpha}\right)\;\Gamma\left(\mu+\frac{r}{\alpha}\right)}{\alpha\,\Gamma(\mu)}\,
      (2^\alpha\delta)^{\frac{r}{\alpha}}+
      \sum_{k=1}^{\infty}\frac{(-1)^k}{k!}\,\frac{\Gamma(\mu+k)}{\Gamma(\mu)}\,\zeta\left(-\alpha
      k+r+1, 1/2\right)\,(2^\alpha\delta)^k\right]\;,
      \end{split}
      \end{equation}
      \begin{equation}\label{assg2}
\begin{split}
&\sum_{k=0}^{\infty}\frac{1-((2k+1)^\alpha\delta+1)^{-\mu}}{(2k+1)^{r+1}}
\underset{\delta\to+0}{\sim}
  \\&
  -2^{-r-1}\left[
  \frac{\Gamma(\mu+p)(-1)^p (2^\alpha\delta)^p}{\Gamma(\mu)\Gamma(p+1)}\left(-\frac{\ln (2^\alpha\delta)}{\alpha}+\frac{\Gamma'(p+1)}{\alpha\,\Gamma(p+1)}-\frac{\Gamma'( 1/2)}{\Gamma(1/2)}-\frac{\Gamma'(\mu+p)}{\alpha\,\Gamma(\mu+p)}\right)
      \right.
      \\&
      \left.+\sum_{k=1,k\ne p}^{\infty}\frac{(-1)^k}{k!}\,\frac{\Gamma(\mu+k)}{\Gamma(\mu)}\,\zeta(-\alpha k+r+1, 1/2)\,(2^\alpha\delta)^k
      \right]
\end{split}
\end{equation}
Пусть $r+1>0$, $\alpha>0$ и $\mu>0$. Тогда имеет место следующее
асимптотическое разложение
\begin{equation}\label{assg3}
\begin{split}
&\sum_{k=0}^{\infty}\frac{1-((2k+1)^\alpha\delta+1)^{-\mu}}{(2k+1)^{r+1}}(-1)^k
\underset{\delta\to+0}{\sim}
\\&
     -2^{-r-1}
      \sum_{k=1}^{\infty}\frac{(-1)^k}{k!}\,\frac{\Gamma(\mu+k)}{\Gamma(\mu)}\,\widetilde{\zeta}(-\alpha k+r+1, 1/2)\,(2^\alpha\delta)^k\;.
      \end{split}
      \end{equation}
\end{example}
\begin{example}\label{SrAriphm}
Пусть $h(t)=(1-t)^{\mu}_+$, $\mu>0$. Второе условие в
теоремах~\ref{thZast1} и \ref{thZast} очевидно выполнено для любых
$\alpha>0$ и $\delta>0$. Если $\mu\ge 1$, то  $m(h)=1$ (см.
замечание~\ref{re_m(h)}). Очевидно $h\in M_m$, $m\in\N$ $\iff\mu\ge
m$. Учитывая замечание~\ref{re_M_m}, получаем следующие результаты.
Если $\alpha\in(0,1]$ и $\mu\ge 2$, то равенство~\eqref{Zast1}
выполняется при любых $n\in\N$, $r\ge\alpha$ и $\delta>0$, а
равенство~\eqref{Zast1_2} выполняется при $n=1$ и
 любых $r\ge\alpha+1$, $\delta>0$.
 Если $\alpha\in(0,1]$ и $\mu\ge 3$, то
равенство~\eqref{Zast1_2} выполняется при любых $n\in\N$,
$r\ge\alpha$ и $\delta>0$. Если $\alpha>1$ и $\mu\ge 1$, то
равенство~\eqref{Zast1} выполняется при  $n=1$ и любых
$r\ge\alpha+1$ и $\delta>0$, а равенство~\eqref{Zast1_2} выполняется
при  $n=1$ и любых $r\ge\alpha+2$ и $\delta>0$. Соответствующие
операторы $R^{\alpha,\mu}_{\delta}=G$ называются средними Рисса. При
$\alpha=\mu=1$, ${1}/{\delta}\in\N$ получаются средние
арифметические.

Второе условие в теореме~\ref{thZast3}  выполнено для любых
$\alpha>0$, $\delta>0$ и $\gamma\in\R$, а значения
$\gamma_m(\rho,h)$ найдены в  работе~\cite{Zast2009_IPMM} (см.
замечание~\ref{re_gamma_m}).
 Поэтому теорема~\ref{thZast3}  справедлива, например,  в следующих
 случаях:
 {\bf 1)}~$\rho=1$,
 $\gamma_1(1,h)={(\mu-2)}/{3}$ (если $\mu\ge 2$),
  $\gamma_2(1,h)={(\mu-3)}/{4}$ (если $\mu\ge 3$).
  {\bf 2)}~$\rho=2$,
 $\gamma_1(2,h)=\mu$ (если $\mu\ge 2$),
  $\gamma_2(2,h)=\mu$ (если $\mu\ge 3$).

Функция $(1-|t|^\alpha)_+^{\mu}$, $\alpha>0$, $\mu>0$, является
положительно определенной на $\R$ $\iff$ $0<\alpha< 2$ и
$\mu\ge\lambda(\alpha)$, где $\lambda(\alpha)$ - функция Кутнера,
которая положительна, возрастает на $(0,2)$, $\lambda(+0)>0$,
$\lambda(1)=1$ и $\lambda(2-0)=+\infty$ (более подробно см.,
например, \cite{Zast2000}). Поэтому операторы
$R^{\alpha,\mu,\gamma}_{\delta}=G$ при $\gamma=0$ будут
положительными при любых $0<\alpha< 2$, $\mu\ge\lambda(\alpha)$ и
$\delta>0$. И наоборот, если при некоторых $\alpha>0$, $\mu>0$
операторы $R^{\alpha,\mu,0}_{\delta}$ будут положительными при любых
$\delta>0$, то $0<\alpha< 2$, $\mu\ge\lambda(\alpha)$. Известно
также, что функция $(1+\gamma|x|)(1-|x|)_+$  является положительно
определенной на $\R$ $\iff$ $\gamma\in[-3,0]$ (см.
\cite[Теорема~9]{Zast2002}). Поэтому положительно определенной на
$\R$ будет и функция $(1+\gamma|x|)(1-|x|)_+^{\mu+1}$ при любых
$\gamma\in[-3,0]$ и $\mu\ge 1$.
 В силу теоремы~\ref{thZast4} равенство~\eqref{Zast2_1} справедливо
 в следующих  случаях:
 {\bf 1)}~$n=1$, $r=1$ или $r\ge 2$, $\gamma=0$, $0<\alpha< 2$, $\mu\ge\lambda(\alpha)$.
 {\bf 2)}~$n=1$, $r=1$ или $r\ge 2$, $\gamma\in[-3,0]$, $\alpha=1$, $\mu=1$ или $\mu\ge 2$.
  Равенство~\eqref{Zast2_2} справедливо
 в следующих  случаях:
 {\bf 1)}  $n=1$, $r=2$ или $r\ge 3$, $\gamma=0$, $0<\alpha< 2$, $\mu\ge\lambda(\alpha)$.
 {\bf 2)}  $n=1$, $r=2$ или $r\ge 3$, $\gamma\in[-3,0]$, $\alpha=1$, $\mu=1$ или $\mu\ge 2$.

 В случае  средних арифметических ($\alpha=\mu=1$, $1/\delta\in\N$) для классов $W^r_p$, $r\in\N$, и
$\widetilde{W}^r_p$, $r-1\in\N$, для которых параметр $n=1$,
равенства~\eqref{Zast1} и ~\eqref{Zast1_2} при $p=\infty$ доказаны в
работах Надя~\cite{Nagy1942,Nagy1946}, а  в работе
Теляковского~\cite{Tel_82} доказано совпадение величин приближения
 указанных классов при $p=1$ и $p=\infty$.
  В этих же случаях асимптотическое разложение при $\delta\to+0$ соответствующих
рядов~\eqref{Zast1} и ~\eqref{Zast1_2} найдено в работах
Теляковского и Баскакова \cite{Tel,Bask-Tel}.
\end{example}

\subsubsection{Приближение средними Чезаро}\label{par_Chesaro}
Числа Чезаро $A_n^{\alpha}$, $n\in\Z_+$, порядка $\alpha\in\R$,
определяются с помощью следующей производящей функции
\begin{equation}\label{ch1}
 \frac{1}{(1-x)^{\alpha+1}}=\sum_{n=0}^{\infty}A_n^\alpha
 x^n\,,\,|x|<1\,.
\end{equation}
Очевидно $A_0^\alpha=1$, $A_k^0=1$ при $k\in\Z_+$, и
$A_n^\alpha=(\alpha+1)\cdot\ldots\cdot(\alpha+n)/n!$  при $n\in\N$.
Естественно считать, что $A_n^\alpha=0$  при $-n\in\N$. Для любых
$\alpha,\gamma\in\R$ и $n\in\Z_+$ справедливы равенства (см.,
например, \cite[гл.~III, \S~1]{Zig})
\begin{equation}\label{ch2}
 A_n^\alpha=\sum_{k=0}^{n}A_k^{\alpha-1}\;\;;\;\;
 A_n^\alpha=\sum_{k=0}^{n}A_{n-k}^{\alpha-\gamma}A_{k}^{\gamma-1}\;.
\end{equation}

 Для заданной последовательности $s_n$, $n\in\Z_+$, чезаровские
 суммы $S_n^\alpha$ порядка $\alpha$ и чезаровские средние $\sigma_n^\alpha$ порядка $\alpha>-1$ определяются по
 формулам
 \begin{equation}\label{ch3}
 S_n^\alpha=\sum_{k=0}^{n}A_{n-k}^{\alpha-1}\,s_k\,,\;\sigma_n^\alpha=\frac{S_n^\alpha}{A_n^\alpha}\,,\,n\in\Z_+\,.
 \end{equation}
 Очевидно $S_n^0=s_n$ и $S_n^1=s_0+\ldots+s_n$, $n\in\Z_+$.
  Нетрудно показать, что для любых
$\alpha,\gamma\in\R$  справедливы равенства (см., например,
\cite[гл.~III, \S~1]{Zig})
\begin{equation}\label{ch4}
 S_n^\alpha=\sum_{k=0}^{n}A_{n-k}^{\alpha-\gamma-1}\,S_{k}^{\gamma}\,,\,n\in\Z_+\,.
\end{equation}
Если в качестве исходной последовательности  взять
 $s_n=D_n(x):=\sum_{\nu=-n}^{n}e^{i\nu x}$ -- ядра Дирихле, то
 чезаровские суммы для этой последовательности будут равны
 \begin{equation}\label{ch5}
 S_n^\alpha(x)=\sum_{k=0}^{n}A_{n-k}^{\alpha-1}D_k(x)=
 \sum_{\nu=-n}^{n}e^{i\nu x}\sum_{k=|\nu|}^{n}A_{n-k}^{\alpha-1}=
 \sum_{\nu=-n}^{n}A_{n-|\nu|}^{\alpha}e^{i\nu x}\,,\,n\in\Z_+\,.
 \end{equation}
 Формула~\eqref{ch4} в этом случае будет иметь вид
 \begin{equation}\label{ch55}
  S_n^\alpha(x)=\sum_{k=0}^{n}A_{n-k}^{\alpha-\gamma-1}\,S_{k}^{\gamma}(x)\,,\,n\in\Z_+\,,\,\alpha,\gamma,x\in\R\,.
 \end{equation}

 Средние Чезаро порядка $\alpha>-1$ функции $f\in L_1$ определяются
 по формуле
 \begin{equation}\label{ch7}
 \begin{split}
 \sigma_n^\alpha(f)(x)=&\sum_{k=-n}^{n}\frac{A_{n-|k|}^\alpha}{A_{n}^\alpha}\widehat{f}(k)e^{ikx}=
 \frac{1}{2\pi}\int_{-\pi}^{\pi}f(t)K_n^\alpha(x-t)\,dt\;;
 \\
 K_n^\alpha(x)=&\frac{1}{A_{n}^\alpha}\cdot
 S_n^\alpha(x)\,,\,n\in\Z_+\,.
 \end{split}
 \end{equation}
 \begin{theorem}\label{thChezaro}
 Пусть $\beta\in\Z$, $p=1$ или $p=\infty$, $\alpha\ge 1$ и $m\in\Z_+$.
  Тогда равенство
  \begin{equation}\label{ch8}
 E(\W^{r,\beta}_{p,1};\sigma_m^\alpha)_p= \displaystyle
 \frac{4}{\pi}\sum\limits_{k=0}^{\infty}\frac{(-1)^{k(\beta+1)}}{(2k+1)^{r+1}}
 \left(1-\frac{A_{m-(2k+1)}^\alpha}{A_m^\alpha}\right)
  \end{equation}
  справедливо по крайней мере в следующих случаях:
  {\bf 1)} $\beta+1\in 2\Z$, $r=1$ или $r\ge 2$;
  {\bf 2)}~$\beta\in 2\Z$, $r=2$ или $r\ge 3$.
  Кроме того, в указанных двух случаях для любых $\alpha,\gamma\ge 1$ и $m\in\Z_+$
   справедливы равенства
   \begin{equation}\label{ch9}
  E(\W^{r,\beta}_{p,1};\sigma_m^\alpha)_p=\frac{1}{A_{m}^{\alpha}}
  \sum_{k=0}^{m}A_{m-k}^{\alpha-\gamma-1}A_{k}^{\gamma}\,
  E(\W^{r,\beta}_{p,1};\sigma_k^\gamma)_p\,.
   \end{equation}
 \end{theorem}
 {\sc Доказательство.}
 Так как
  \begin{equation*}
  S_m^1(x)=D_0(x)+\ldots+D_m(x)=\frac{\sin^2\frac{(m+1)x}{2}}{\sin^2\frac{x}{2}}\ge 0
  \;,\;m\in\Z_+\;,\;x\in\R\,,
  \end{equation*}
  то из~\eqref{ch55} при $\gamma=1$ вытекает, что любых $\alpha>1$,
  $m\in\Z_+$ и
  $x\in\R$ выполняется неравенство $S_m^\alpha(x)>0$. Ряд Фурье для
  ядра $K\in L_1$ (см.~\eqref{K}) в нашем случае имеет вид
  \begin{equation*}
K(t)\sim
  \sum_{k\ne 0}\frac{\left(1-\frac{A_{m-|k|}^\alpha}{A_m^\alpha}\right)e^{-i\beta\pi\sign  k/2}}{|k|^r}\,e^{ikt}=
  \sum_{k=1}^{\infty}\frac{2\left(1-\frac{A_{m-k}^\alpha}{A_m^\alpha}\right)}{k^r}\,
  \cos\left(kt-\frac{\beta\pi}{2}\right)\;.
\end{equation*}
  В случае $\beta+1\in 2\Z$ надо применить пример~\ref{ex1}~({\bf iii}) и   равенство~\eqref{fs1},
 а в  случае $\beta\in 2\Z$ надо применить пример~\ref{ex2}~({\bf iv}) и  равенство~\eqref{fc1},
  в которых
 $\nu_k=A_{m-k}^\alpha/A_{m}^\alpha$, $k\in\Z_+$ и $|c|=2$.
 При этом надо учесть, что $\nu_0=1$ и
 \begin{equation*}
 \frac{\nu_0}{2}+\sum_{k=1}^{\infty}\nu_k\cos
 kt=\frac{1}{2A_{m}^\alpha}\cdot S_m^\alpha(t)\ge
 0\,,\,m\in\Z_+\,,\,\alpha\ge 1\,,\,t\in\R\,.
 \end{equation*}
 Равенство~\eqref{ch8} в указанных двух случаях доказано.
 Докажем теперь, что в этих же случаях справедливы
 равенства~\eqref{ch9}.

Для фиксированных параметров  $\beta\in\Z$ и
  $r> 0$ зададим две последовательности
 $s_m=s_m(\beta,r)$ и $E_m^\alpha=E_m^\alpha(\beta,r)$, $m\in\Z_+$, по следующему правилу:
 \begin{equation*}
 {s}_{2k}=0\;,\;{s}_{2k+1}=\frac{4}{\pi}\cdot\frac{(-1)^{k(\beta+1)}}{(2k+1)^{r+1}}\,,\,
 k\in\Z_+\;;\;
 E_m^{\alpha}=S A_m^\alpha-S_m^{\alpha+1}\,,\,m\in\Z_+\,,
 \end{equation*}
 где $S_m^\alpha$ -- чезаровские
 суммы последовательности $s_m$, $m\in\Z_+$, а
  $S=\sum_{k=0}^{\infty}s_k$. Из~\eqref{ch4} сразу получаются
  равенства
\begin{equation}\label{ch6}
 \sum_{k=0}^{m}A_{m-k}^{\alpha-\gamma-1}E_k^\gamma=E_m^\alpha
 \;,\;\alpha,\gamma\in\R\,,\,m\in\Z_+\,.
   \end{equation}
  Кроме того, из~\eqref{ch3} вытекает, что
  \begin{equation*}
  E_m^\alpha=
 \sum_{k=0}^{\infty}s_k(A_m^\alpha-A_{m-k}^\alpha)
  =\displaystyle
 \frac{4}{\pi}\sum\limits_{k=0}^{\infty}\frac{(-1)^{k(\beta+1)}}{(2k+1)^{r+1}}
 \left({A_m^\alpha}-{A_{m-(2k+1)}^\alpha}\right)\,,\,m\in\Z_+\,,\,\alpha\in\R\,.
 \end{equation*}
 Если $p=1$ или $p=\infty$, а параметры $\beta\in\Z$ и
  $r>0$ такие, для которых выполняется одно из двух указанных в
  теореме  условий, то по доказанному при любых $\alpha\ge 1$ и $m\in\Z_+$
  справедливо равенство~\eqref{ch8}, т.е.
   $E_m^\alpha=A_m^\alpha E(\W^{r,\beta}_{p,1};\sigma_m^\alpha)_p$.
   В этих случаях равенство~\eqref{ch9} вытекает из
   равенства~\eqref{ch6}.
  Теорема~\ref{thChezaro} доказана.
  \begin{remark}
  Для классов $W^r_p$, $r\in\N$, и
$\widetilde{W}^r_p$, $r-1\in\N$, равенства~\eqref{ch8} для
$\alpha\in\N$ при $p=\infty$ доказаны в работах
Надя~\cite{Nagy1942,Nagy1946}, а в работе Теляковского~\cite{Tel_82}
доказано совпадение  при $p=1$ и $p=\infty$ величин приближения
 указанных классов средними Чезаро.
  \end{remark}
   \begin{remark}
   Если $\beta\in\Z$, а $p=1$ или $p=\infty$, то для средних Фейера  равенство
 \begin{equation}\label{ch10}
 E(\W^{r,\beta}_{p,1};\sigma_m^1)_p=
 \frac{K^{r,\beta}}{m+1}
 +O\left(\frac{1}{(m+1)^r}\right)\;,\;\text{где},\;
 K^{r,\beta}=\frac{4}{\pi}\sum\limits_{k=0}^{\infty}\frac{(-1)^{k(\beta+1)}}{(2k+1)^{r}}\,,
 \end{equation}
   справедливо
по крайней мере в следующих случаях:
  {\bf 1)} $\beta+1\in 2\Z$,  $r\ge 2$;
  {\bf 2)}~$\beta\in 2\Z$, $r=2$ или $r\ge 3$.
 Для $r-1\in\N$ равенство~\eqref{ch10} доказано Никольским~\cite{Nik1941,Nik1945}.
 Для остальных указанных $r$ равенство~\eqref{ch10} получается из
 равенства~\eqref{ch8} методом, изложенным  в работе
Теляковского и Баскакова \cite{Bask-Tel}. Из равенства~\eqref{ch9}
при $\gamma=1$ вытекает, что при  $\alpha>1$ и указанных выше
значениях $p$, $\beta$, $r$  справедливо равенство
 \begin{equation*}
  E(\W^{r,\beta}_{p,1};\sigma_m^\alpha)_p=\frac{1}{A_{m}^{\alpha}}
  \sum_{k=0}^{m}A_{m-k}^{\alpha-2}(k+1)
  \left(\displaystyle
  \frac{K^{r,\beta}}{k+1}+\frac{\varepsilon_k}{k+1}\right)=
  \frac{(K^{r,\beta}+\widetilde{\sigma}_m^{\alpha-1})\alpha}{\alpha+m}\,,
   \end{equation*}
  где
   $\varepsilon_k=  (k+1)E(\W^{r,\beta}_{p,1};\sigma_k^1)_p-
   {K^{r,\beta}}$,
   а $\widetilde{\sigma}_m^{\gamma}=
  \sum_{k=0}^{m}A_{m-k}^{\gamma-1}\varepsilon_k/{A_{m}^{\gamma}}$ -- средние Чезаро
  порядка $\gamma>-1$ для последовательности $\varepsilon_k$, $k\in\Z_+$.
  Так как $\widetilde{\sigma}_m^{0}=\varepsilon_m\to 0$, то
  (см., например, \cite[гл.~III, \S~1]{Zig})
  $\widetilde{\sigma}_m^{\alpha-1}\to 0$ при $\alpha>1$ и, значит,
  $E(\W^{r,\beta}_{p,1};\sigma_m^\alpha)_p={K^{r,\beta}\alpha}/{(\alpha+m)}+o(1/{(\alpha+m)})$.
  Последнее соотношение при $p=\infty$ и $r-1\in\N$ другим методом
  доказано Фалалеевым~\cite{Falaleev2007}.
  \end{remark}

\subsubsection{Приближение средними типа Рисса и Чезаро целого порядка}\label{par_Polinom}
Пусть $Q(x)=\sum_{k=1}^{n}a_kx^{\mu_k}$, где $n\in\N$ и
$a_k,\mu_k>0$. Так как $Q(0)=0$ и $Q(x)>0$ при $x>0$, то для любых
$\alpha,u>0$ имеет смысл следующий оператор
\begin{equation}\label{Q1}
 \begin{split}
 G_{u}^{\alpha,Q}(f)(x):=&\sum_{k\in\Z}
  \frac{Q\left((u-|k|^\alpha)_+\right)}{Q(u)}\,\widehat{f}(k)e^{ikx}=
  \frac{1}{2\pi}\int_{-\pi}^{\pi}f(t)K_u^{\alpha,Q}(x-t)\,dt\,,\,f\in
  L_1\,,
  \\
 K_u^{\alpha,Q}(x)=&\sum_{k\in\Z}\frac{Q\left((u-|k|^\alpha)_+\right)}{Q(u)}\,e^{ikx}\,.
 \end{split}
\end{equation}
 Если $Q(x)=x^\mu$, $\mu>0$, то получаем
 $G_{u}^{\alpha,Q}=R^{\alpha,\mu}_{1/u}$ -- средние Рисса (см. пример~\ref{SrAriphm}).
  Если $Q(x)=\prod_{k=1}^{\mu}(x+k-1)$, $\mu\in\N$, то при $\alpha=1$, $m\in\Z_+$, получаем
 $G_{m+1}^{1,Q}=\sigma^{\mu}_{m}$ -- средние Чезаро (см. \S~\ref{par_Chesaro}).

 Для той же функции $Q$  и для любых
$\alpha,\delta>0$ определим следующие операторы
 \begin{equation}
 \begin{split}
 \widetilde{G}_{\delta}^{\alpha,Q}(f)(x):=&\sum_{k\in\Z}
  \frac{Q\left((1-|k|^\alpha\delta)_+\right)}{Q(1)}\,\widehat{f}(k)e^{ikx}=
  \frac{1}{2\pi}\int_{-\pi}^{\pi}f(t)\widetilde{K}_{\delta}^{\alpha,Q}(x-t)\,dt\,,\,f\in
  L_1\,,
  \\
 \widetilde{K}_{\delta}^{\alpha,Q}(x)=&\sum_{k\in\Z}\frac{Q\left((1-|k|^\alpha\delta)_+\right)}{Q(1)}\,e^{ikx}\,.
 \end{split}
\end{equation}
  \begin{theorem}\label{thPolinom}
Пусть $Q(x)=\sum_{k=1}^{n}a_kx^{\mu_k}$, где $n\in\N$,  $a_k>0$ и
$\mu_k\ge 1$.
  Пусть $\beta\in\Z$, $p=1$ или $p=\infty$, $0<\alpha\le 1$ и $u,\delta>0$.
  Тогда равенства
  \begin{equation}\label{Q3}
 E(\W^{r,\beta}_{p,1};G_{u}^{\alpha,Q})_p= \displaystyle
 \frac{4}{\pi}\sum\limits_{k=0}^{\infty}\frac{(-1)^{k(\beta+1)}}{(2k+1)^{r+1}}
 \left(1-\frac{Q\left((u-(2k+1)^\alpha)_+\right)}{Q(u)}\right)\,,
  \end{equation}
    \begin{equation}\label{Q4}
 E(\W^{r,\beta}_{p,1};\widetilde{G}_{\delta}^{\alpha,Q})_p= \displaystyle
 \frac{4}{\pi}\sum\limits_{k=0}^{\infty}\frac{(-1)^{k(\beta+1)}}{(2k+1)^{r+1}}
 \left(1-\frac{Q\left((1-(2k+1)^\alpha\delta)_+\right)}{Q(1)}\right)
  \end{equation}
  справедливы по крайней мере в следующих случаях:
  {\bf 1)} $\beta+1\in 2\Z$, $r=1$ или $r\ge 2$;
  {\bf 2)}~$\beta\in 2\Z$, $r=2$ или $r\ge 3$.
  \end{theorem}
   {\sc Доказательство.}
   Рассмотрим случай операторов $G_{u}^{\alpha,Q}$.
   При любых $u>0$,
   $0<\alpha\le 1$, $\mu\ge 1$, функция $g(t)=(u-|t|^\alpha)^\mu_+$  является положительно определенной на $\R$ (более подробно см.,
например, \cite{Zast2000}).
   Поэтому при любых $u>0$,
    $0<\alpha\le 1$, положительно определенной на $\R$ будет и функция
    $f(t)=Q((u-|t|^\alpha)_+)$. По лемме~\ref{le0} при всех $t\in\R$
    выполняется неравенство $ K_u^{\alpha,Q}(t)\ge 0$.
 Ряд Фурье для
  ядра $K\in L_1$ (см.~\eqref{K}) в нашем случае имеет вид
  \begin{equation*}
K(t)\sim
  \sum_{k\ne 0}\frac{\left(1-\frac{Q\left((u-|k|^\alpha)_+\right)}{Q(u)}\right)e^{-i\beta\pi\sign  k/2}}{|k|^r}\,e^{ikt}=
  \sum_{k=1}^{\infty}\frac{2\left(1-\frac{Q\left((u-|k|^\alpha)_+\right)}{Q(u)}\right)}{k^r}\,
  \cos\left(kt-\frac{\beta\pi}{2}\right)\;.
\end{equation*}
    В случае $\beta+1\in 2\Z$ надо применить пример~\ref{ex1}~({\bf iii}) и   равенство~\eqref{fs1},
 а в  случае $\beta\in 2\Z$ надо применить пример~\ref{ex2}~({\bf iv}) и  равенство~\eqref{fc1},
  в которых
 $\nu_k={Q\left((u-k^\alpha)_+\right)}/{Q(u)}$, $k\in\Z_+$ и $|c|=2$.
  При этом надо учесть, что $\nu_0=1$ и
 \begin{equation*}
 \frac{\nu_0}{2}+\sum_{k=1}^{\infty}\nu_k\cos
 kt=\frac{1}{2}\cdot K_u^{\alpha,Q}(t)\ge
 0\,,\,u>0\,,\,0<\alpha\le 1\,,\,t\in\R\,.
 \end{equation*}
  Точно так же рассматривается и случай операторов
  $\widetilde{G}_{\delta}^{\alpha,Q}$.
   Теорема~\ref{thPolinom} доказана.

\subsection{Вспомогательные утверждения}\label{par_Vspom}
\subsubsection{Свойства выпуклых функций и величины $m(h)$}
Отметим следующие два свойства выпуклых функций.\\
{\bf 1.} Если обе функции $f$ и $g$ неотрицательны, убывают и
выпуклы вниз на $(0,+\infty)$, то такой же функцией будет и
произведение $fg$ (доказательство вытекает из определения выпуклой
функции и очевидного неравенства $(f(x)-f(y))(g(x)-g(y))\ge 0$ при
всех $x,y>0$).\\
{\bf 2.} Если функция $f$ неотрицательна, убывает и выпукла вниз на
$(0,+\infty)$, а функция $g$ положительна, возрастает и выпукла
вверх на $(0,+\infty)$, то  функция $f(g(x))$ неотрицательна,
убывает и выпукла вниз на $(0,+\infty)$ (доказательство легко
получается из определения  выпуклой функции). В частности, при любом
$\varepsilon\in(0,1]$ функция $f(x^\varepsilon)$ неотрицательна,
убывает и выпукла вниз на $(0,+\infty)$.
\begin{lemma}\label{le3}
Пусть функция $\lambda(x)$ выпукла вниз на $(0,+\infty)$ и
$\lambda(t_k)\to 0$ для некоторой последовательности
$t_k\to+\infty$. Тогда при любых $\varepsilon\in(0,1]$,
 $\gamma\ge 0$
функция $\lambda(x^\varepsilon)x^{-\gamma}$ убывает к нулю и выпукла
вниз на $(0,+\infty)$.
 Если дополнительно $\lambda(x)\in C^1(0,+\infty)$ и
функция $-\lambda'(x)$ выпукла вниз на $(0,+\infty)$, то при любых
$\varepsilon\in(0,1]$, $\gamma\ge 0$
 функция
$-f'(x)$ монотонно убывает к нулю и выпукла вниз на $(0,+\infty)$,
где $f(x)=\lambda(x^\varepsilon)x^{-\gamma}$.
\end{lemma}
{\sc Доказательство.} Пусть функция $\lambda$ выпукла вниз на
$(0,+\infty)$. Тогда справедливо неравенство
\begin{equation}\label{vypukl}
\frac{\lambda(x_2)-\lambda(x_1)}{x_2-x_1}\le
\frac{\lambda(x_3)-\lambda(x_1)}{x_3-x_1}\le
\frac{\lambda(x_3)-\lambda(x_2)}{x_3-x_2}\;,\;0<x_1<x_2<x_3\;.
\end{equation}
Если в неравенстве~\eqref{vypukl} взять $x_3=t_k$ и перейти к
пределу при $k\to\infty$, то получим, что функция $\lambda(x)$
убывает на $(0,+\infty)$  и, значит,
$\lim\limits_{x\to+\infty}\lambda(x)=\lim\limits_{k\to+\infty}\lambda(t_k)=0$.
   Поэтому при любых $\varepsilon\in(0,1]$, $\gamma\ge 0$
    функция $\lambda(x^\varepsilon)x^{-\gamma}$ убывает к нулю  и выпукла вниз на
    $(0,+\infty)$.

Пусть дополнительно $\lambda(x)\in C^1(0,+\infty)$ и функция
$-\lambda'(x)$ выпукла вниз на интервале $(0,+\infty)$. Из равенства
$\lambda(+\infty)=0$ вытекает сходимость интеграла
$\int_{1}^{+\infty}\lambda'(x)\,dx$. Но функция $\lambda'(x)$
монотонна на интервале $(x_0,+\infty)$ при некотором $x_0\ge 0$.
Поэтому $\lambda'(+\infty)=0$ и, значит, функция $-\lambda'(x)$
убывает к нулю на $(0,+\infty)$. Тогда при любых
$\varepsilon\in(0,1]$, $\gamma\ge 0$
 функция
$-f'(x)=\varepsilon(-\lambda'(x^\varepsilon))
x^{-(\gamma+1-\varepsilon)}+\gamma\lambda(x^\varepsilon)
x^{-\gamma-1}$
 монотонно убывает к нулю и выпукла вниз на $(0,+\infty)$.
 Лемма~\ref{le3} доказана.

В следующей лемме установлены свойства  величины $m(h)$ (см.
определение~\ref{def_m(h)}).
\begin{lemma}\label{le2}
Пусть $h(x)\le 1$ при  $x>0$. Тогда справедливы следующие
утверждения:\\
{\bf 1.}
  Если   $m(h)<+\infty$ и $h(x_0)=1$ при некотором $x_0>0$, то $h(x)=1$ при всех $x\ge x_0$.
  \\
  {\bf 2.}
  $m(h)=-\infty\iff h(x)\equiv 1$ на $(0,+\infty)$.
\\
{\bf 3.} Если $-\infty<m(h)<+\infty$, то функция
$(1-h(x^\varepsilon))x^{-\gamma}$, $\varepsilon>0$, убывает по
$x\in(0,+\infty)$
 $\iff \gamma\ge\varepsilon m(h)$.
  \\
  {\bf 4.}
 Если $h(+0)=1$ и $h(x)\not\equiv 1$ на $(0,+\infty)$, то $m(h)>0$.
 \\
 {\bf 5.}
 Если $h(x_k)\to q<1$ для некоторой последовательности $x_k\to+\infty$, то $m(h)\ge 0$.
 Если, дополнительно, $h(x_0)>q$ в некоторой точке $x_0>0$,  то $m(h)>0$.
 \\
 {\bf 6.}
 Если $h'$ непрерывна в точке $0$ справа, $h(0)=1$ и $h'(0)\ne 0$, то $m(h)\ge 1$.
\\
{\bf 7.}
 Если функция $h$ выпукла вниз на $(0,+\infty)$ и $h(x)\not\equiv const$ на $(0,+\infty)$,
 то $m(h)\in(0,1]$.
 \\
 {\bf 8.}
  Пусть $x_0>0$ такая точка, что $h(x)<1$ при $0<x<x_0$
   и $h(x)=1$ при $x\ge x_0$, а  если $h(x)<1$ при всех $x>0$, то считаем $x_0:=+\infty$.
   Если функция $h$ дифференцируема на $(0,x_0)$,
   то $m(h)=\sup\{h'(x)\,x/(h(x)-1):0<x<x_0\}$.
\end{lemma}
{\sc Доказательство.} При некотором $\gamma\in\R$ функция
$f(x)=(1-h(x))x^{-\gamma}$ убывает и неотрицательна на
$(0,+\infty)$. Поэтому, если $f(x_0)=0$ при некотором $x_0>0$, то
$f(x)=0$ при всех $x\ge x_0$. Утверждение 1 доказано.

 Докажем утверждение 2. Пусть  $m(h)=-\infty$. Тогда при всех $\gamma\in\R$ функция
$f(x)=(1-h(x))x^{-\gamma}$ убывает  на $(0,+\infty)$. Поэтому
функция $f$ а, значит, и функция $h$ дифференцируемы при почти всех
$x\in(0,+\infty)$ (относительно меры Лебега). В тех точках
$x\in(0,+\infty)$, где дифференцируема функция $h$ следующее
неравенство будет выполняться при всех  $\gamma\in\R$:
\begin{equation}\label{proizvod}
f'(x)=-x^{-\gamma-1}\left(xh'(x)+\gamma(1-h(x))\right)\le 0\iff
xh'(x)+\gamma(1-h(x))\ge 0\;.
\end{equation}
Если обе части последнего неравенства разделить на $|\gamma|\ne 0$ и
перейти к пределу при $\gamma\to-\infty$, то получим неравенство
$h(x)\ge 1$ и, значит, $h(x)= 1$ во всех точках  $x\in(0,+\infty)$,
где дифференцируема функция $h$. Так как множество таких точек всюду
плотно на $(0,+\infty)$, то учитывая утверждение 1, получаем, что
$h(x)\equiv 1$ на $(0,+\infty)$. Необходимость в утверждении~2
доказана. Достаточность очевидна.

Докажем утверждение 3. Если $-\infty<m(h)<+\infty$, то очевидно
функция $(1-h(x^\varepsilon))x^{-\gamma}$, $\varepsilon>0$, убывает
по $x\in(0,+\infty)$
 $\iff$ функция $(1-h(x))x^{-\gamma/\varepsilon}$ убывает по $x\in(0,+\infty)$
  $\iff \gamma/\varepsilon\ge m(h)$.

Докажем утверждение 4. Если $m(h)\le 0$, то функция $1-h(x)$ убывает
на $(0,+\infty)$. Поэтому $0\le 1-h(x)\le 1-h(+0)=0$ при все $x>0$,
что противоречит условию $h(x)\not\equiv 1$ на $(0,+\infty)$.

Докажем утверждение 5. Если при некотором $\gamma<0$ функция
$\lambda(x)=(1-h(x))x^{-\gamma}$ убывает  на $(0,+\infty)$, то
существует конечный предел $\lambda(+\infty)$, но
$\lambda(x_k)\to+\infty$. Поэтому $m(h)\ge 0$. Если функция $1-h(x)$
убывает  на $(0,+\infty)$, то функция $h(x)$ возрастает и, значит,
$h(x)\le h(+\infty)=q$ для всех $x>0$. Поэтому, если дополнительно
$h(x_0)>q$ в некоторой точке $x_0>0$,  то $m(h)>0$. Утверждение 5
доказано.

Докажем утверждение 6. Пусть при некотором $\gamma\in\R$ функция
$f(x)=(1-h(x))x^{-\gamma}$ убывает на $(0,+\infty)$. Из условия
вытекает, что при некотором $x_0>0$ неравенство~\eqref{proizvod}
выполняется при всех $x\in(0,x_0)$. Неравенство~\eqref{proizvod}
делим на $x>0$ и переходим к пределу при $x\to+0$. Получим
неравенство $h'(0)(1-\gamma)\ge 0$. Так как $h'(0)\ne 0$, то
$h'(0)<0$ (иначе $h(x)>h(0)=1$ при малых  $x>0$, что противоречит
условию). Поэтому $\gamma\ge 1$ и, значит, $m(h)\ge 1$. Утверждение
6 доказано.

Докажем утверждение 7. Пусть функция $h$ выпукла вниз на
$(0,+\infty)$. Тогда справедливо неравенство
\begin{equation}\label{vypukl_2}
\frac{h(x_2)-h(x_1)}{x_2-x_1}\le \frac{h(x_3)-h(x_1)}{x_3-x_1}\le
\frac{h(x_3)-h(x_2)}{x_3-x_2}\;,\;0<x_1<x_2<x_3\;.
\end{equation}
Так как функция $h$ ограничена сверху, то из
неравенства~\eqref{vypukl_2} вытекает, что функция $h$ убывает на
$(0,+\infty)$ и, значит, существует конечный предел $h(+0)\le 1$.
Если в неравенстве~\eqref{vypukl} перейти к пределу при $x_1\to+0$,
то получим, что функция $(h(+0)-h(x))/x$  убывает на $(0,+\infty)$.
Следовательно убывает на $(0,+\infty)$ и функция $(1-h(x))/x$.
 Поэтому
$m(h)\le 1$. Если $m(h)\le 0$, то функция  $1-h(x)$ убывает  на
$(0,+\infty)$, т.е. $h(x)$ возрастает и, значит, $h(x)\equiv const$
на $(0,+\infty)$, что противоречит условию. Поэтому $m(h)>0$.
Утверждение 7 доказано.

Докажем утверждение 8. При любом $\gamma\in\R$ функция
$f(x)=(1-h(x))x^{-\gamma}$ неотрицательна на $(0,+\infty)$ и
$f(x)=0$ при $x\ge x_0$ (если $x_0<+\infty$). Поэтому функция $f$
убывает  на $(0,+\infty)\iff$  $f$ убывает  на $(0,x_0)\iff$ при
всех $x\in(0,x_0)$ выполняется неравенство~\eqref{proizvod} $\iff$
$\gamma\ge \sup\{h'(x)\,x/(h(x)-1):0<x<x_0\}$. Утверждение 8
доказано. Лемма~\ref{le2} полностью доказана.

 \subsubsection{Классы $M_m$ и свойства величины $\gamma_m(\rho,h)$}\label{subMm}
\begin{theorem}\label{Williamson}
Пусть $m\in\N$. Тогда следующие три условия эквивалентны:\\
{\bf 1.}
 $h\in C^{m-1}(0,+\infty)$ и  функции $(-1)^kh^{(k)}$,
$k=0,1,\dots,m-1$, неотрицательны, убывают и выпуклы вниз на
$(0,+\infty)$.
\\
{\bf 2.} $h\in C^{m-1}(0,+\infty)$,  функция $(-1)^{m-1}h^{(m-1)}$
неотрицательна, убывает, выпукла вниз на $(0,+\infty)$ и существует
конечный предел $h(+\infty)\ge 0$.
\\
{\bf 3.} Существует неотрицательная борелевская мера $\mu$ на
$[0,+\infty)$, для которой величина $\mu([0,a])$ конечна для всех
$a>0$ и
\begin{equation}\label{Schonberg}
h(x)=\int_{0}^{+\infty}(1-xs)^m_+\,d\mu(s)\;,\;x>0\;.
\end{equation}
\end{theorem}
Из представления~\ref{Schonberg} следует, что
$h(+\infty)=\mu(\{0\})$ и $h(+0)=\mu([0,+\infty))$. Функции,
удовлетворяющие условию~1 в теореме~\ref{Williamson}, называются
кратно монотонными. Формулу \eqref{Schonberg} для кратно монотонных
функций в 1940~г. получил Schoenberg. Доказательство
теоремы~\ref{Williamson} можно найти в работе Williamson~\cite{Will}
(см. также работу L\'evy~\cite{Levy}).

 В определении~\ref{defM_m} класс функций, удовлетворяющих
условию~2 в теореме~\ref{Williamson} мы обозначили через $M_m$,
$m\in\N$. Из теоремы~\ref{Williamson}, вытекает, что $M_{m+1}\subset
M_{m}$. Отметим следующие два свойства (см.~\cite{Will}):
 {\bf 1)} Если $f,g\in M_m$, то и $fg\in M_m$.
 {\bf 2)} Если $f\in M_m$, то $\lim\limits_{x\to+\infty}x^k
 f^{(k)}(x)=0$ при $k=1,\ldots,m$ (у выпуклой функции $(-1)^{m-1}f^{(m-1)}(x)$
  во всех точках существуют односторонние производные, которые совпадают почти всюду).

Введем следующие функции
\begin{equation*}
\varphi_m(x):=(m+1)\int_{0}^{1}(1-xs)^m_+\,ds
 \,,\,x\ge 0\,,\;m\in\N\;.
\end{equation*}
Из теоремы~\ref{Williamson} вытекает, что $\varphi_m\in M_m$.
Очевидно $\varphi_m\in C^{m-1}[0,+\infty)$ и
\begin{equation}
\varphi_m(x)=\frac{1-(1-x)^{m+1}_+}{x}\,,\,x>0\,;\,
 \varphi^{(k-1)}_m(x)=\frac{(-1)^{k-1}(k-1)!}{x^k}\,,\,x>1\,,\,k\in\N\,.
\end{equation}
Поэтому для любой неотрицательной конечной борелевской меры $\mu$ на
$[0,+\infty)$
\begin{equation}\label{M_m}
f(x):=\int_{0}^{+\infty}\varphi_m(xs)s\,d\mu(s)\in M_m\,.
\end{equation}
\begin{sledstvie}\label{sl}
Если $h\in M_{m+1}$ при некотором  $m\in\N$ и существует конечный
предел $h(+0)\le 1$, то $\lambda(x)=(1-h(x))/x\in M_k$ при всех
$k=1,\ldots,m$.
\end{sledstvie}
{\sc Доказательство.} Для функции $h$ имеет место
представление~\eqref{Schonberg} (с заменой $m$ на $m+1$). В нашем
случае мера $\mu$ из этого представления конечна и
$h(+0)=\mu([0,+\infty))$. Учитывая~\eqref{M_m}, получаем, что
\begin{equation*}
 \lambda(x)=\frac{1-h(+0)}{x}+\frac{h(+0)-h(x)}{x}=
 \frac{1-h(+0)}{x}+\int_{0}^{+\infty}\varphi_{m}(xs)s\,d\mu(s)\in M_{m}\,.
\end{equation*}
Осталось учесть, что  $M_{m}\subset M_{k}$ при всех $k=1,\ldots,m$.
\begin{remark}
 Очевидно $f\in M:=\bigcap_{m\in\N}M_m$ $\iff$ $f\in C^{\infty}(0,+\infty)$ и неравенство $f^{(k)}(x)\ge 0$
 выполняется для всех $k\in\Z_+$ и $x>0$.
 Такие функции называются вполне монотонными на $(0,+\infty)$.
 Теорема Бернштейна-Хаусдорфа-Уиддера~\cite{Bernstein,Hausdorff,Widder} утверждает, что функция $f$
 вполне монотонна на $(0,+\infty)$ ($f\in M$) $\iff$
 $f(x)=\int_0^{+\infty}e^{-xt}\ d\mu(t)$, $ x>0$, где $\mu$ неотрицательная
  борелевская мера на $[0,+\infty)$ такая, что интеграл сходится
  для всех $x>0$. Из следствия~\ref{sl} вытекает, что функция
  $\lambda(x)=(1-h(x))/x\in M$ для любой функции $h\in M$ с
  условием $h(+0)\le1$.
 \end{remark}
В следующей лемме установлены свойства  величины $\gamma_m(\rho,h)$
и функции $\lambda_{\rho,\gamma}(x)$ (см.
определение~\ref{def_gamma_m}).
\begin{lemma}\label{le5}
Если при некотором  $m\in\N$ функция $h\in M_{m+1}$, $0<h(+0)\le 1$
и $h(+\infty)=0$, то справедливы следующие утверждения:\\
{\bf 1.}
 Если при некоторых $\rho\ge 1$, $\gamma\in\R$ функция
$(-1)^{m-1}\lambda^{(m-1)}_{\rho,\gamma}(x)$ выпукла вниз на
$(0,+\infty)$, то $(-1)^{m-1}\lambda^{(m-1)}_{\rho,\gamma}(x)$
неотрицательна и убывает к нулю на $(0,+\infty)$ и, значит, функция
$\lambda_{\rho,\gamma}(x)\in M_{m}$.
\\
{\bf 2.}
Если $\rho\ge 1$, то $0\le\gamma_m(\rho,h)<+\infty$.
\\
{\bf 3.}
 Если $\rho\ge 1$, то
$(-1)^{m-1}\lambda^{(m-1)}_{\rho,\gamma}(x)$ выпукла вниз на
$(0,+\infty)$ $\iff$ $\lambda_{\rho,\gamma}(x)\in M_{m}$ $\iff
\gamma\le\gamma_m(\rho,h)$.
\\
{\bf 4.}
 Функция $\gamma_m(\rho,h)$ возрастает по
$\rho\in[1,+\infty)$.
\\
{\bf 5.}
 Если дополнительно $h\in M_{m+2}$, то
$\gamma_{m+1}(\rho,h)\le\gamma_m(\rho,h)$ при всех $\rho\ge 1$.
\end{lemma}
{\sc Доказательство.}
 Из следствия~\ref{sl} вытекает, что обе функции
$(1-h(x))/x$ и $h(x)$ принадлежат классу $M_m$.
 Так как $x^{-\varepsilon}\in M_m$ при  всех $\varepsilon\ge 0$,
 то  обе функции
$f(x)=(1-h(x))/x^\rho$ и $g(x)=h(x)/x^{\rho-1}$ принадлежат классу
$M_m$ при всех $\rho\ge 1$. Отсюда следует, что $\gamma_m(\rho,h)\ge
0$. Из свойств функций из $M_m$ вытекает, что $f^{(k)}(+\infty)=0$ и
$g^{(k)}(+\infty)=0$ при всех $k=1,\ldots,m$, а если $h(+\infty)=0$,
то и при $k=0$. Поэтому $\lambda^{(m-1)}_{\rho,\gamma}(+\infty)=0$
при любых $\rho\ge 1$ и $\gamma\in\R$. Отсюда следует утверждение 1.
Равенство $\gamma_m(\rho,h)=+\infty$ не возможно, ибо в противном
случае функция $h(x)/x^{\rho-1}$ будет одновременно выпуклой вниз и
вверх на  $(0,+\infty)$, что противоречит условиям $h(+\infty)=0$ и
$h(+0)>0$. Утверждение 3 вытекает из первого. Утверждение 4 вытекает
из того, что произведение двух функций из $M_m$ также является
функцией  из $M_m$. Утверждение~5 вытекает из вложений
$M_{m+2}\subset M_{m+1}\subset M_{m}$. Лемма~\ref{le5} доказана.

\subsection{Доказательство теорем}\label{par_Proof}

\subsubsection{Доказательство  теоремы~\ref{thZast1}}
  Ряд Фурье для ядра $K\in L_1$  (см.~\eqref{K})  в нашем
 случае имеет вид
\begin{equation}
K(t)\sim
  \sum_{k\ne 0}\frac{(1-h(|k|^\alpha\delta))e^{-i\beta\pi\sign  k/2}}{|k|^r}\,e^{ikt}=
  \sum_{k=1}^{\infty}\frac{2(1-h(k^\alpha\delta))}{k^r}\,
  \cos\left(kt-\frac{\beta\pi}{2}\right)\;.
\end{equation}

1) Если $\beta=2p+1$, $p\in\Z$, то
\begin{equation}\label{sin}
  (-1)^p K(t)\sim\sum_{k=1}^{\infty}\frac{2(1-h(k^\alpha\delta))}{k^r}\,\sin kt=
  \sum_{k=1}^{\infty}\lambda_k\sin kt\;.
\end{equation}
 К функции $\lambda(x)$ применяем лемму~\ref{le3} при $t_k=k^\alpha\delta$
 (последовательность $h(|k|^\alpha\delta)$, как коэффициенты Фурье
интегрируемой функции $g_{\alpha,\delta}$, стремится к нулю).
  Так как $\alpha\in(0,1]$, то при
 $r\ge\alpha$   функция $f(x)=2\delta\lambda(x^\alpha\delta)x^{\alpha-r}$
 убывает к нулю и выпукла вниз на $(0,+\infty)$. Поэтому последовательность $\lambda_k=f(k)$, $k\in\N$,
 монотонно убывает к нулю и выпукла вниз.
 В силу примера~\ref{ex1}~({\bf i})  равенство~\eqref{fs} справедливо при всех $n\in\N$.
 Первое утверждение в теореме~\ref{thZast1} доказано.

2) Если  $\beta=2p$, $p\in\Z$, то
\begin{equation}\label{cos}
  (-1)^p K(t)\sim\sum_{k=1}^{\infty}\frac{2(1-h(k^\alpha\delta))}{k^r}\,\cos kt=
  \sum_{k=1}^{\infty}\mu_k\cos kt\;.
\end{equation}
Так как $\alpha\in(0,1]$, то при
 $r\ge\alpha+1$  последовательность $k\mu_k=f(k)$, $k\in\N$,
 где $f(x)=2\delta\lambda(x^\alpha\delta)x^{\alpha+1-r}$,
 монотонно убывает к нулю и выпукла вниз. В силу примера~\ref{ex2}~({\bf iii})  равенство~\eqref{fc} справедливо
при  $n=1$.

 Пусть дополнительно $\lambda(x)\in C^1(0,+\infty)$ и
функция $-\lambda'(x)$ выпукла вниз на интервале $(0,+\infty)$. Так
как $\alpha\in(0,1]$, то при $r\ge\alpha$ последовательность
$\mu_k=f(k)$, $k\in\N$,
  убывает к нулю и выпукла вниз, где $f(x)=2\delta\lambda(x^\alpha\delta)x^{\alpha-r}$
 и функция $-f'(x)$  выпукла вниз на $(0,+\infty)$ (см.  лемму~\ref{le3}). Поэтому
 $\Delta^3\mu_k=\int_{k}^{k+1}(2f'(x+1)-f'(x)-f'(x+2))\,dx\ge 0$ для всех $k\in\N$.
В силу примера~\ref{ex2}~({\bf i})  равенство~\eqref{fc} справедливо
при всех $n\in\N$.  Теорема~\ref{thZast1}  доказана.



\subsubsection{Доказательство  теоремы~\ref{thZast}}
 Неравенства $m(h)\ge 0$ и $m(h)>0$, если $h(x_0)>0$ при некотором $x_0>0$,
 вытекают из утверждения 5 в лемме~\ref{le2} при $x_k=k^\alpha\delta$, $k\in\N$, и $q=0$
 (последовательность $h(|k|^\alpha\delta)$, как коэффициенты Фурье
интегрируемой функции $g_{\alpha,\delta}$, стремится к нулю).

1) Пусть $\beta=2p+1$, $p\in\Z$ и $r\ge\alpha  m(h)+1$. В этом
 случае для ядра $K\in L_1$ имеет место~\eqref{sin}, где
$\lambda_k\ge 0$ и последовательность $k\lambda_k$ убывает (см.
утверждение 3 в лемме~\ref{le2}). В силу примера~\ref{ex1}~({\bf
ii}) равенство~\eqref{fs} справедливо при $n=1$.

2) Пусть  $\beta=2p$, $p\in\Z$ и  $r\ge\alpha  m(h)+2$. В этом
 случае для ядра
$K\in L_1$ имеет место~\eqref{cos}, где $\mu_k\ge 0$ и
последовательность $k^2\mu_k$ убывает (см. утверждение 3 в
лемме~\ref{le2}). В силу примера~\ref{ex2}~({\bf ii})
равенство~\eqref{fc} справедливо при $n=1$. Теорема~\ref{thZast}
 доказана.
\subsubsection{Доказательство теоремы~\ref{thZast3}}
 Ряд Фурье для ядра $K\in L_1$  (см.~\eqref{K})  в нашем
 случае имеет вид
\begin{equation}
K(t)\sim
  \sum_{k=1}^{\infty}\frac{2(1-(1+k^\alpha\delta\gamma)h(k^\alpha\delta))}{k^r}\,
  \cos\left(kt-\frac{\beta\pi}{2}\right)\;.
\end{equation}

1) Если $\beta=2p+1$, $p\in\Z$, то
\begin{equation}\label{sin2}
  (-1)^p K(t)\sim\sum_{k=1}^{\infty}\frac{2(1-(1+k^\alpha\delta\gamma)h(k^\alpha\delta))}{k^r}\,\sin kt=
  \sum_{k=1}^{\infty}\lambda_k\sin kt\;.
\end{equation}
 Здесь $\lambda_k=f(k)$, $k\in\N$, где
  $f(x)=2\delta^\rho\lambda_{\rho,\gamma}(x^\alpha\delta)x^{\alpha\rho-r}$,
  $x>0$, а функция $\lambda_{\rho,\gamma}(x)$ определена равенством~\eqref{lambda_rho}.
 Если $\alpha\in(0,1]$ и
 $r\ge\alpha\rho$, то   функция $f(x)$
 убывает к нулю и выпукла вниз на $(0,+\infty)$. Поэтому последовательность $\lambda_k=f(k)$, $k\in\N$,
 монотонно убывает к нулю и выпукла вниз.
 В силу примера~\ref{ex1}~({\bf i})  равенство~\eqref{fs} справедливо при всех $n\in\N$.
 Если $\alpha>1$ и
 $r\ge\alpha\rho+1$, то   функция $xf(x)$
 убывает к нулю  на $(0,+\infty)$. Поэтому последовательность $\lambda_k=kf(k)$, $k\in\N$,
 монотонно убывает к нулю.
 В силу примера~\ref{ex1}~({\bf ii})  равенство~\eqref{fs} справедливо при $n=1$.
 Первое утверждение в теореме~\ref{thZast3} доказано.

2) Если  $\beta=2p$, $p\in\Z$, то
\begin{equation}\label{cos2}
  (-1)^p K(t)\sim\sum_{k=1}^{\infty}\frac{2(1-(1+k^\alpha\delta\gamma)h(k^\alpha\delta))}{k^r}\,\cos kt=
  \sum_{k=1}^{\infty}\mu_k\cos kt\;.
\end{equation}
Здесь $\mu_k=f(k)$, $k\in\N$, где, как и в первом случае,
  $f(x)=2\delta^\rho\lambda_{\rho,\gamma}(x^\alpha\delta)x^{\alpha\rho-r}$,
  $x>0$.
Если $\alpha\in(0,1]$ и
 $r\ge\alpha\rho+1$, то функция $xf(x)$,
 монотонно убывает к нулю и выпукла вниз. Поэтому последовательность $k\mu_k$,
 $k\in\N$, монотонно убывает к нулю и выпукла вниз.
 В силу примера~\ref{ex2}~({\bf iii})  равенство~\eqref{fc} справедливо при  $n=1$.
 Если $\alpha>1$ и
 $r\ge\alpha\rho+2$, то функция $x^2f(x)$,
 монотонно убывает к нулю. Поэтому последовательность $k^2\mu_k$,
 $k\in\N$, монотонно убывает и неотрицательна.
 В силу примера~\ref{ex2}~({\bf ii})  равенство~\eqref{fc} справедливо при  $n=1$.

 Если $\alpha\in(0,1]$ и
 $r\ge\alpha\rho$, то   функция $f(x)=2\delta^\rho\lambda_{\rho,\gamma}(x^\alpha\delta)x^{\alpha\rho-r}$
 убывает к нулю и выпукла вниз на $(0,+\infty)$ (см.
лемму~\ref{le3}).
 Поэтому последовательность $\mu_k=f(k)$, $k\in\N$,
 монотонно убывает к нулю и выпукла вниз.
 Если дополнительно  $h\in M_3$ и  $\gamma\le\gamma_2(\rho,h)$, то
функция $-\lambda'_{\rho,\gamma}(x)$ выпукла вниз на $(0,+\infty)$
и, значит,    функция $-f'(x)$  выпукла вниз на $(0,+\infty)$ (см.
лемму~\ref{le3}). Поэтому
 $\Delta^3\mu_k=\int_{k}^{k+1}(2f'(x+1)-f'(x)-f'(x+2))\,dx\ge 0$ для всех $k\in\N$.
В силу примера~\ref{ex2}~({\bf i})  равенство~\eqref{fc} справедливо
при всех $n\in\N$.
 Теорема~\ref{thZast3}  доказана.
\subsubsection{ Доказательство теоремы~\ref{thZast4}}
    Функция
 $f(x):=(1+\gamma |t|^\alpha\delta)h(|t|^\alpha\delta)$
 является  положительно определенной и непрерывной на $\R$.
   Мы рассматриваем случай~\eqref{sin2}, где
   $\lambda_k=2(1-\nu_k)/k$, $k\in\N$, и $\nu_k=f(k)$, $k\in\Z_+$.
   В силу леммы~\ref{le0} и примера~\ref{ex1}~({\bf iii})  равенство~\eqref{fs1} справедливо при $n=1$, $r=1$ и $r\ge 2$.
   В случае~\eqref{cos2}, где
   $\mu_k=2(1-\nu_k)/k$, $k\in\N$, и $\nu_k=f(k)$, $k\in\Z_+$ из примера~\ref{ex2}~({\bf iv}) получаем,
   что  равенство~\eqref{fc1} справедливо при $n=1$, $r=2$ и $r\ge 3$.
   Теорема~\ref{thZast4} доказана.


\end{document}